# A Lagrangian approach for solving an axisymmetric thermo-electromagnetic problem. Application to time-varying geometry processes

Marta Benítez[1,2†], Alfredo Bermúdez[1,3†], Pedro Fontán[4†], Iván Martínez[1,3†], Pilar Salgado[1,3*†]

[1]CITMAga, Galician Centre for Mathematical Research and Technology, Santiago de Compostela, E-15782, Spain.
[2]Department of Mathematics, Universidade da Coruña, Elviña s/n, A Coruña, 15071, Spain.
[3]Departament of Applied Mathematics, Universidade de Santiago de Compostela, Santiago de Compostela, E-15782, Spain.
[4]Repsol Technology Lab, Spain.

*Corresponding author(s). E-mail(s): mpilar.salgado@usc.es;
Contributing authors: marta.benitez@udc.es; alfredo.bermudez@usc.es; pedro.fontan@repsol.com; ivanmartinez.suarez@usc.es;
†These authors contributed equally to this work.

**Abstract**

The aim of this work is to introduce a thermo-electromagnetic model for calculating the temperature and the power dissipated in cylindrical pieces whose geometry varies with time and undergoes large deformations; the motion will be a known data. The work will be a first step towards building a complete thermo-electromagnetic-mechanical model suitable for simulating electrically assisted forming processes, which is the main motivation of the work. The electromagnetic model will be obtained from the time-harmonic eddy current problem with an in-plane current; the source will be given in terms of currents or voltages defined at some parts of the boundary. Finite element methods based on a Lagrangian weak formulation will be used for the numerical solution. This approach will avoid the need to compute and remesh the thermo-electromagnetic domain along the time. The numerical tools will be implemented in FEniCS and validated by using a suitable test also solved in Eulerian coordinates.

# 1 Introduction

The aim of this work is to describe mathematical models and numerical tools to simulate the thermo-electromagetic behaviour of cylindrical ferromagnetic pieces whose geometry can vary over time. The study is motivated by a particular type of Electrically Assisted Forming process for preforming cylindrical bars that undergo large deformations in short time.

Electrically Assisted Forming (EAF) is a processing technique which applies electricity during the plastic deformation of materials. This technology has gained large interest in the industry for manufacturing pieces in automotive or aeroplane sectors. We are interested in one of these techniques, known as electric upsetting, which is a preforming process able to create a local enlarged diameter at one end of a bar to be later forged without further heating. In this process, a cold bar is placed in an horizontal upsetter machine and clamped by gripper jaws. The electric current passes through the bar by contact between one of its ends and the gripper jaws. The bar end is heated up and it acquires a plastic behaviour. When it reaches enough temperature, the bar is pushed against the anvil with the help of a force applied by a pusher located at the opposite end and the diameter at the hot end is enlarged. Fig. 1 shows the evolution of the piece during the upsetting and Fig. 2 shows a simple sketch and the main elements involved in the process.

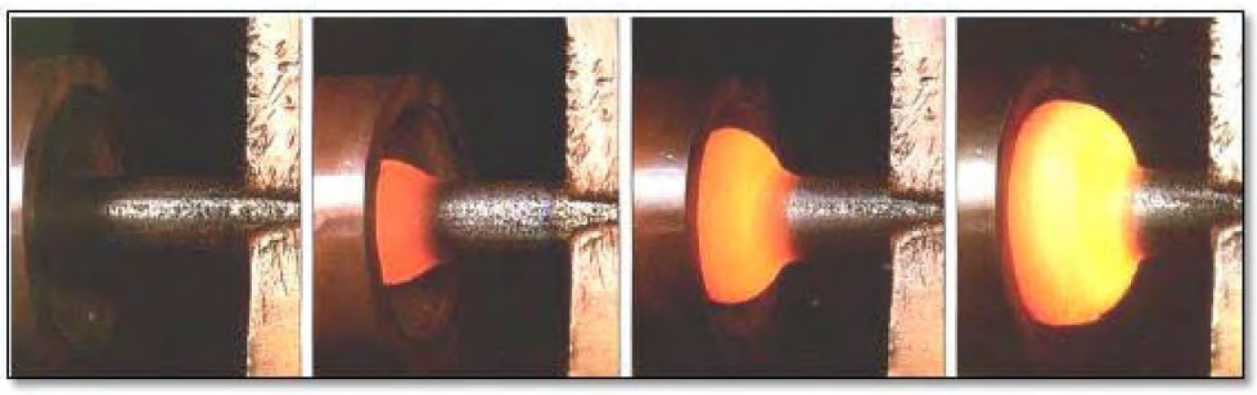

**Fig. 1**: Enlarging the diameter at the end of a steel bar submitted to electric upsetting.

Deformation and heating occur simultaneously and it is essential to control both the electrical input signal and the upsetting force at the same time to achieve the desired shape and quality [1]. The full simulation involves a thermo-electromagnetic-mechanical model and the proposals in the literature [1, 2] are usually based on a sequential multiphysics approach where the electromagnetic model is restricted to a steady linear case due to the consideration of a direct current source. However, if the industrial process is powered by alternating current and the material is ferromagnetic, a non-linear eddy current model is required. Our aim is to construct a fully thermo-electromagnetic-mechanical model allowing the consideration of eddy currents and non-linear magnetic materials. To achieve this goal, in this paper we make a first step focusing on building a thermo-electromagnetic model which assumes that the deformation is known. We will adopt a Lagrangian approach, which avoids computing and remeshing the computational domain along time. Although this approach is widely used in structural mechanics, it is more unusual in coupled electromagnetic-mechanical systems with large deformations. We refer the reader to [3] for a Lagrangian approach to medical EAF problems; see also in [4] a work focused on electromagnetic forming which deals with a Lagrangian approach by using a least-action variational principle. A detailed description of Maxwell's equations in material form can be found in [5] (see also [6]). It is worth mentioning that the Lagrangian formalism has been exploited by Bossavit to compute forces in deformables bodies [7, 8]. On the other hand, the authors of the present paper have recently made the first steps to consider the Lagrangian approach to solve thermo-electro-mechanical problems by using a direct current source and assuming small deformations [9].

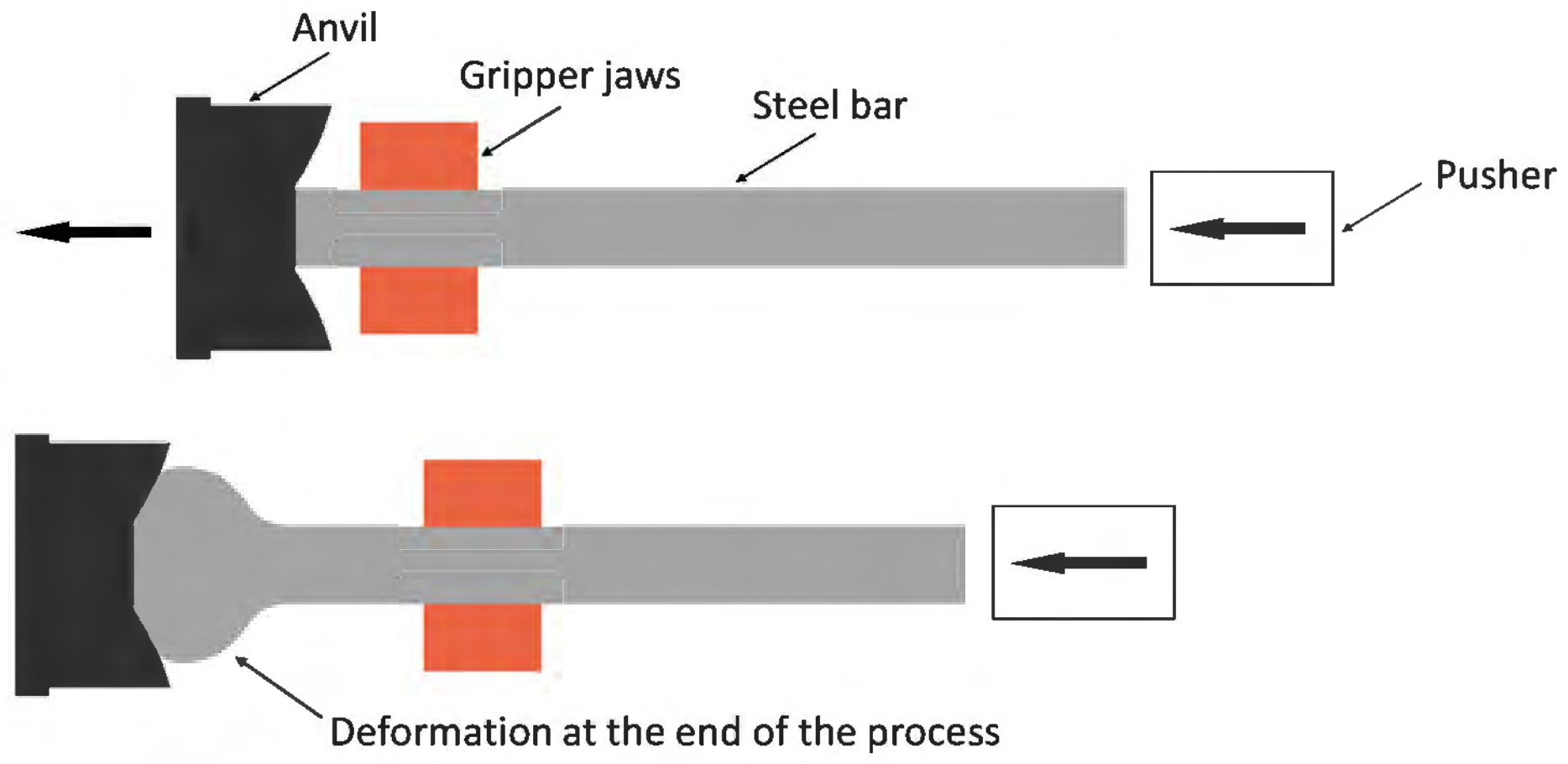

**Fig. 2**: Sketch of the electric upsetting process.

The thermo-electromagnetic problem that arises when modelling the electric upsetting process requires to solve an eddy current model to determine the power dissipated at the different parts of the domain that are heated up and a transient heat transfer model to determine the temperature along the time. As mentioned earlier, in this paper we will assume that the deformation is a data of the problem. The coupling between the thermal and the electromagnetic problems is due to the temperature dependence of the electromagnetic properties and to the heat source provided by the solution of the electromagnetic model. We will study the coupled model in an axisymmetric framework, considering that the current source is usually given in terms of either the total current at some electrical ports of the electro upsetter which are connected to a transformer, or the potential drops between such ports. The eddy current model with this kind of boundary conditions is referred to in the literature as the *eddy current model with electrical ports* or *with non-local boundary conditions* and has been extensively studied in 3D domains. A pioneering work is due to A. Bossavit [10] which was the starting point for many authors who later studied the problem using different formulations ([6, 11]) often based on the ideas introduced by Bossavit [12–15], consisting in the use of different unknown fields in the eddy current model (electric field, magnetic field or combining these with scalar potentials).

On the other hand, the axisymmetric eddy current model with electrical ports and with an in-plane current has recently been studied from a mathematical and numerical point of view in a linear framework [16]. An interesting feature of this model is that it can be written in terms of the azimuthal component of the magnetic field and can be formulated in the conducting domain by defining suitable boundary conditions without the need to consider the surrounding air. Now, with a time-dependent domain, the Eulerian formulation of the eddy current model will be presented and then the Lagrangian formulation will be derived using tools from continuum mechanics. To attain this goal, the 3D formulation will be first established and then written in the axisymetric case. In this paper, materials with non-linear and temperature-dependent magnetic laws will also be considered.

The paper is structured as follows: Section 2 introduces some notations and definitions concerning the motion of a continuum body, which are needed for the later description of the Lagrangian approaches; this section also introduces the main tools for describing the models in an axisymmetric setting. In Section 3 we describe the electromagnetic model for a bounded domain with cylindrical symmetry and electrical ports, first in Eulerian coordinates and then in the reference configuration. In Section 4 we follow the same steps to introduce the thermal model. Finally, in Section 5 we describe the discretization tools used to solve the coupled problem and present an academic test to validate the numerical tools; in particular, we solve the thermo-electromagnetic model in a cylinder with a known displacement field that emulates the real deformation of an electric upsetting process.

## 2 Preliminary tools and notation

Let $\Omega$ be a bounded domain in $\mathbb{R}^3$ with Lipschitz boundary $\Gamma$. Let $\mathsf{X} : \bar{\Omega} \times \mathbb{R} \to \mathbb{R}^3$ be a *motion* in the sense of Gurtin [17]. In particular, $\mathsf{X} \in C^3(\bar{\Omega} \times \mathbb{R})$ and, for each fixed $t \in \mathbb{R}$, $\mathsf{X}(\cdot, t)$ is a one-to-one function satisfying

$$\det \mathbf{F} > 0, \qquad \text{in } \bar{\Omega} \times \mathbb{R}, \tag{1}$$

being $\mathbf{F}(:, t) = \mathbf{Grad}\, \mathsf{X}(:, t)$ the deformation gradient tensor. Notice that $\bar{\Omega}(t) := \mathsf{X}(\bar{\Omega}, t)$ is a closed region for all $t$. In practice, a bounded time interval is considered, namely, $[0, T]$.

For a material point $\mathsf{p} \in \Omega$, its *position* at time $t$ is given by $\mathsf{x} = \mathsf{X}(\mathsf{p}, t)$ and its *velocity* in the material configuration is given by the time derivative of the motion, that is $\dot{\mathsf{X}}(\mathsf{p}, t)$.

Let us introduce the *trajectory* of the motion

$$\mathcal{T} := \{(\mathsf{x}, t) : \mathsf{x} \in \bar{\Omega}(t),\ t \in [0, T]\}. \tag{2}$$

We remark that fields defined in $\mathcal{T}$ (respectively, in $\bar{\Omega} \times [0, T]$) are called spatial fields (respectively, material fields). If $\phi$ is a spatial field we define its material

description $\phi_m$ by

$$\phi_m(\mathsf{p}, t) := \phi(\mathsf{X}(\mathsf{p}, t), t). \tag{3}$$

The spatial description of the velocity is $\mathbf{v}(\mathsf{x}, t) = \dot{\mathsf{X}}(\mathsf{P}(\mathsf{x}, t), t)$, being $\mathsf{P} : \mathcal{T} \to \bar{\Omega}$ the so-called *reference map* of the motion, that is the inverse of the one-to-one mapping $\mathsf{X}(:, t)$, defined by

$$\mathsf{p} = \mathsf{P}(\mathsf{x}, t) \Leftrightarrow \mathsf{x} = \mathsf{X}(\mathsf{p}, t). \tag{4}$$

We can assume, without loss of generality, that for $t = 0$ $\mathsf{X}(\mathsf{p}, 0) = \mathsf{p}$. Thus, $\mathsf{p} = \mathsf{P}(\mathsf{x}, t)$ can be considered as the initial position of the point that occupies the place $\mathsf{x}$ at time $t$.

By following the same ideas, we can define the material *acceleration*, $\ddot{\mathsf{X}}(\mathsf{p}, t)$, and the spatial description of the acceleration, $\mathbf{a}(\mathsf{x}, t) = \ddot{\mathsf{X}}(\mathsf{P}(\mathsf{x}, t), t)$.

Let us denote by $\mathbf{u}(\mathsf{p}, t)$ the *material displacement* of point $\mathsf{p}$ at time $t$, namely,

$$\mathbf{u}(\mathsf{p}, t) = \mathsf{X}(\mathsf{p}, t) - \mathsf{p}. \tag{5}$$

On the other hand, for a smooth scalar spatial field $\phi(\mathsf{x}, t)$, the *material time derivative* is given by

$$\dot{\phi}(\mathsf{x}, t) = \frac{\partial \phi}{\partial t}(\mathsf{x}, t) + \mathbf{v}(\mathsf{x}, t) \cdot \mathbf{grad}\phi(\mathsf{x}, t). \tag{6}$$

Table 1 collects the notation used for the Eulerian (spatial) and Lagrangian (material) descriptions in the paper.

**Table 1**: Notation for the Eulerian and Lagrangian descriptions.

| **Element** | **Eulerian** | **Lagrangian** |
|---|---|---|
| Domain | $\Omega(t)$ | $\Omega$ |
| Domain boundary | $\Gamma(t)$ | $\Gamma$ |
| Normal unit vector | $\mathbf{n}_{\mathsf{x}}$ | $\mathbf{n}_{\mathsf{p}}$ |
| Tangent unit vector | $\boldsymbol{\tau}_{\mathsf{x}}$ | $\boldsymbol{\tau}_{\mathsf{p}}$ |
| Differential line element | $\mathrm{d}l_{\mathsf{x}}$ | $\mathrm{d}l_{\mathsf{p}}$ |
| Differential surface element | $\mathrm{d}A_{\mathsf{x}}$ | $\mathrm{d}A_{\mathsf{p}}$ |
| Differential volume element | $\mathrm{d}V_{\mathsf{x}}$ | $\mathrm{d}V_{\mathsf{p}}$ |
| Magnitude | $L(\mathsf{x}, t)$ | $L_m(\mathsf{p}, t)$ |
| Differential gradient operator | **grad** | **Grad** |
| Differential curl operator | **curl** | **Curl** |

We will also need the following identities to transform integrals defined in the Eulerian description to the Lagrangian one. Let $\phi(\mathsf{x}, t)$ be a scalar spatial field and $\boldsymbol{\vartheta}(\mathsf{x}, t)$ a spatial vector field. Then

$$\int_{\Gamma(t)} \phi(\mathsf{x}, t) \mathrm{d}A_{\mathsf{x}} = \int \phi_m(\mathsf{p}, t) \det \mathbf{F}(\mathsf{p}, t) |\mathbf{F}^{-t}(\mathsf{p}, t) \mathbf{n}_{\mathsf{p}}| \, \mathrm{d}A_{\mathsf{p}}, \tag{7}$$

$$\int_{\Omega(t)} \phi(\mathsf{x}, t) \mathrm{d}V_{\mathsf{x}} = \int_{\Omega} \phi_m(\mathsf{p}, t) \det \mathbf{F}(\mathsf{p}, t) \, \mathrm{d}V_{\mathsf{p}}, \tag{8}$$

$$\int_{(t)} \boldsymbol{\vartheta}(\mathsf{x},t) \cdot \mathbf{n}_{\mathsf{x}} \mathrm{d}A_{\mathsf{x}} = \int_{\Gamma} \boldsymbol{\vartheta}_m(\mathsf{p},t) \cdot \det \mathbf{F}(\mathsf{p},t) \mathbf{F}^{-t}(\mathsf{p},t) \mathbf{n}_{\mathsf{p}} \, \mathrm{d}A_{\mathsf{p}}, \tag{9}$$

where $\mathbf{n}_{\mathsf{p}}$ (respect. $\mathbf{n}_{\mathsf{x}}$) is an outward unit normal vector to the boundary of $\Omega$ at point $\mathsf{p}$ (respect. to the boundary of $\Omega(t)$ at point $\mathsf{x}$).

Next, we introduce some notation related to problems with cylindrical symmetry. Let us assume that the reference domain $\Omega$ can be obtained by rotating a bounded domain $\widehat{\Omega}$, with boundary $\partial\widehat{\Omega} = \widehat{\Gamma}_{\mathrm{D}} \cup \widehat{\Gamma}$, around the axis of symmetry:

$$\Omega := \{(r_m, \theta_m, z_m) : \theta_m \in [0, 2\pi), (r_m, z_m) \in \widehat{\Omega}\}, \tag{10}$$

and

$$\Gamma := \partial\Omega = \{(r_m, \theta_m, z_m) : \theta_m \in [0, 2\pi), (r_m, z_m) \in \widehat{\Gamma}\}. \tag{11}$$

Notice that $\widehat{\Gamma}_{\mathrm{D}}$ is defined by

$$\widehat{\Gamma}_{\mathrm{D}} := \{(r_m, z_m) \in \partial\widehat{\Omega} : r_m = 0\}. \tag{12}$$

Let us consider $\phi_m(\mathsf{p},t)$, a material scalar field, $\boldsymbol{\vartheta}_m(\mathsf{p},t)$, a material vector field, and $\boldsymbol{\Psi}_m(\mathsf{p},t)$, a material tensor field, all of them with cylindrical symmetry. Notice that $\mathsf{p}$ is a point of the reference configuration $\Omega$ with radial and axial coordinates given by $\widehat{\mathsf{p}} := (r_m, z_m) \in \widehat{\Omega}$. Let us define

$$\widehat{\phi}_m(\widehat{\mathsf{p}},t) := \phi_m(\mathsf{p},t), \tag{13}$$

and $\widehat{\boldsymbol{\vartheta}}_m(\widehat{\mathsf{p}},t)$ and $\widehat{\boldsymbol{\Psi}}_m(\widehat{\mathsf{p}},t)$ defined by their coordinates in the basis $\{\mathbf{e}_r, \mathbf{e}_\theta, \mathbf{e}_z\}$ as follows

$$\widehat{\boldsymbol{\vartheta}}_m(\widehat{\mathsf{p}},t) = \begin{pmatrix} \widehat{\vartheta}_r(\widehat{\mathsf{p}},t) \\ 0 \\ \widehat{\vartheta}_z(\widehat{\mathsf{p}},t) \end{pmatrix} := \boldsymbol{\vartheta}_m(\mathsf{p},t), \tag{14}$$

$$\widehat{\boldsymbol{\Psi}}_m(\widehat{\mathsf{p}},t) = \begin{pmatrix} \widehat{\Psi}_{rr}(\widehat{\mathsf{p}},t) & 0 & \widehat{\Psi}_{rz}(\widehat{\mathsf{p}},t) \\ 0 & \widehat{\Psi}_{\theta\theta}(\widehat{\mathsf{p}},t) & 0 \\ \widehat{\Psi}_{zr}(\widehat{\mathsf{p}},t) & 0 & \widehat{\Psi}_{zz}(\widehat{\mathsf{p}},t) \end{pmatrix} := \boldsymbol{\Psi}_m(\mathsf{p},t). \tag{15}$$

Moreover, let us introduce the material tensor field $\widehat{\widehat{\boldsymbol{\Psi}}}_m(\widehat{\mathsf{p}},t)$ defined by

$$\widehat{\widehat{\boldsymbol{\Psi}}}_m(\widehat{\mathsf{p}},t) = \begin{pmatrix} \widehat{\Psi}_{rr}(\widehat{\mathsf{p}},t) & \widehat{\Psi}_{rz}(\widehat{\mathsf{p}},t) \\ \widehat{\Psi}_{zr}(\widehat{\mathsf{p}},t) & \widehat{\Psi}_{zz}(\widehat{\mathsf{p}},t) \end{pmatrix}. \tag{16}$$

Notice that the previous relations can be defined in a similar way for spatial fields with cylindrical symmetry.

Therefore, if $\widehat{\mathbf{u}}(\widehat{\mathsf{p}},t) = \widehat{u}_r(\widehat{\mathsf{p}},t)\widehat{\mathbf{e}}_r + \widehat{u}_z(\widehat{\mathsf{p}},t)\widehat{\mathbf{e}}_z$ is the displacement field associated to point $\widehat{\mathsf{p}}$ at time $t$, the position of $\widehat{\mathsf{p}}$ in the current configuration (i.e., in the Eulerian configuration) will be $\widehat{\mathsf{x}} = (r,z)$ with components

$$r = r_m + \widehat{u}_r(\widehat{\mathsf{p}},t), \tag{17}$$

$$z = z_m + \widehat{u}_z(\widehat{\mathsf{p}},t). \tag{18}$$

For a problem with cylindrical symmetry, the deformation gradient tensor $\mathbf{F}(\mathsf{p},t)$ has the form

$$\mathbf{F}(\mathsf{p},t) = \mathbf{I_3} + \mathbf{Grad}\,\mathbf{u}(\mathsf{p},t) := \widehat{\mathbf{F}}(\widehat{\mathsf{p}},t) = \begin{pmatrix} 1+\dfrac{\partial \widehat{u}_r(\widehat{\mathsf{p}},t)}{\partial r_m} & 0 & \dfrac{\partial \widehat{u}_r(\widehat{\mathsf{p}},t)}{\partial z_m} \\ 0 & 1+\dfrac{\widehat{u}_r(\widehat{\mathsf{p}},t)}{r_m} & 0 \\ \dfrac{\partial \widehat{u}_z(\widehat{\mathsf{p}},t)}{\partial r_m} & 0 & 1+\dfrac{\partial \widehat{u}_z(\widehat{\mathsf{p}},t)}{\partial z_m} \end{pmatrix}.$$

Let us introduce the material differential operator $\widehat{\mathbf{Grad}}$ which applied to a scalar material field with cylindrical symmetry $\widehat{\phi}_m(\widehat{\mathsf{p}},t)$ is defined by

$$\widehat{\mathbf{Grad}}\,\widehat{\phi}_m(\widehat{\mathsf{p}},t) = \begin{pmatrix} \dfrac{\partial \widehat{\phi}_m(\widehat{\mathsf{p}},t)}{\partial r_m} \\ \dfrac{\partial \widehat{\phi}_m(\widehat{\mathsf{p}},t)}{\partial z_m} \end{pmatrix}, \tag{19}$$

and, applied to a vector field with cylindrical symmetry, $\widehat{\boldsymbol{\vartheta}}_m(\widehat{\mathsf{p}},t)$, provides

$$\widehat{\mathbf{Grad}}\,\widehat{\boldsymbol{\vartheta}}_m(\widehat{\mathsf{p}},t) = \begin{pmatrix} \dfrac{\partial \widehat{\vartheta}_r(\widehat{\mathsf{p}},t)}{\partial r_m} & \dfrac{\partial \widehat{\vartheta}_r(\widehat{\mathsf{p}},t)}{\partial z_m} \\ \dfrac{\partial \widehat{\vartheta}_z(\widehat{\mathsf{p}},t)}{\partial r_m} & \dfrac{\partial \widehat{\vartheta}_z(\widehat{\mathsf{p}},t)}{\partial z_m} \end{pmatrix}. \tag{20}$$

In a similar way, it is possible to define the spatial operator $\widehat{\mathbf{grad}}$.

The definition (20) leads to the introduction of the tensor $\widehat{\widehat{\mathbf{F}}}(\widehat{\mathsf{p}},t)$, defined by

$$\widehat{\widehat{\mathbf{F}}}(\widehat{\mathsf{p}},t) = \mathbf{I_2} + \widehat{\mathbf{Grad}}\,\widehat{\mathbf{u}}(\widehat{\mathsf{p}},t) = \begin{pmatrix} 1+\dfrac{\partial \widehat{u}_r(\widehat{\mathsf{p}},t)}{\partial r_m} & \dfrac{\partial \widehat{u}_r(\widehat{\mathsf{p}},t)}{\partial z_m} \\ \dfrac{\partial \widehat{u}_z(\widehat{\mathsf{p}},t)}{\partial r_m} & 1+\dfrac{\partial \widehat{u}_z(\widehat{\mathsf{p}},t)}{\partial z_m} \end{pmatrix}, \tag{21}$$

which verifies the relation

$$\det \mathbf{F}(\mathsf{p},t) = \left(1+\frac{\widehat{u}_r(\widehat{\mathsf{p}},t)}{r_m}\right)\det \widehat{\widehat{\mathbf{F}}}(\widehat{\mathsf{p}},t). \tag{22}$$

# 3 The electromagnetic model

In this section we describe the electromagnetic model which allows us to compute the electromagnetic fields in a piece which is deformed during a time interval, assuming that the deformation is known. The electrical source will be provided in terms of voltages and currents through some parts of the boundary. In particular, this model will be suitable for computing the electromagnetic fields in a cylindrical piece which is deformed by using an electric upsetting machine operating with alternating current. We will introduce, in a first step, a 3D model with electric ports to obtain later the axisymmetric model with its corresponding Lagrangian formulation. An interesting feature is that the cylindrical symmetry allows us to state the problem only in the conducting part, which will be denoted by $\Omega(t)$, by using suitable boundary conditions.

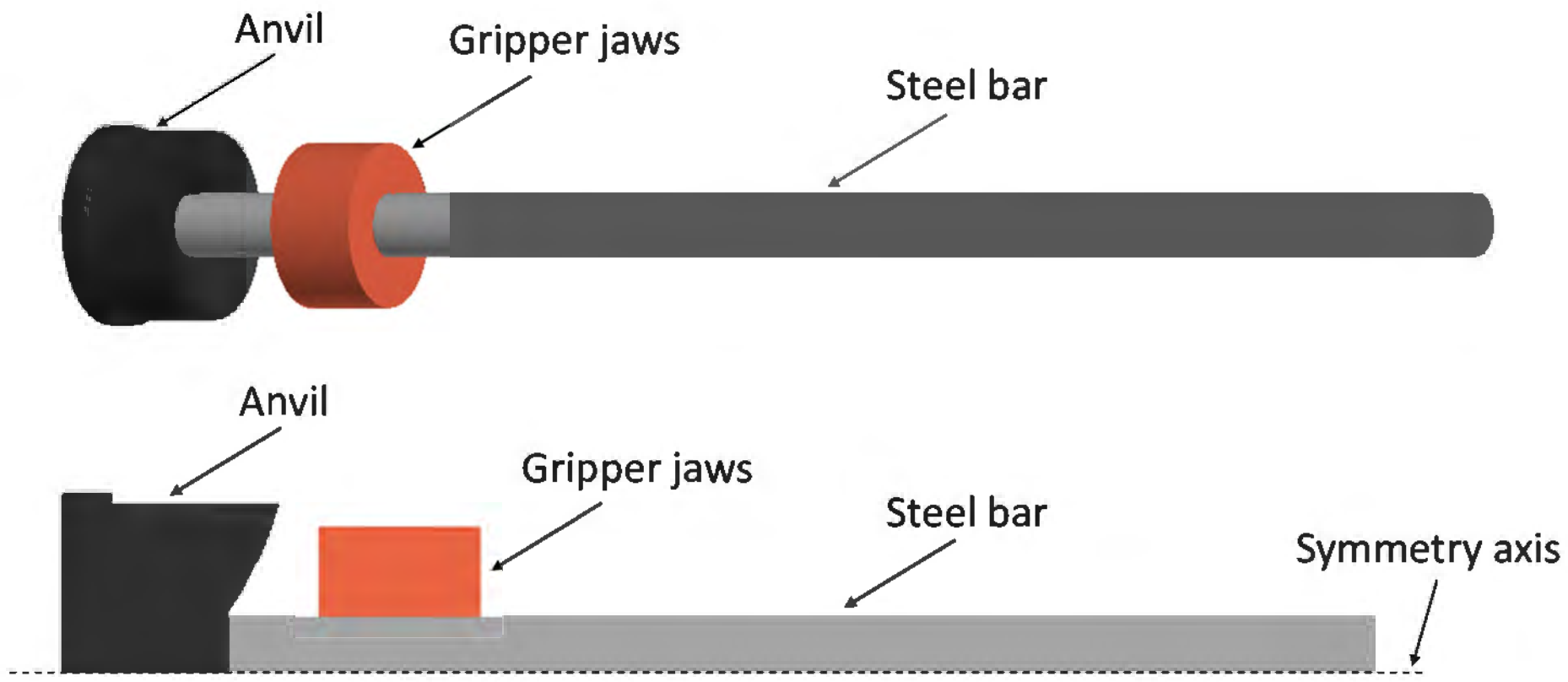

**Fig. 3**: Sketch of the elements considered in the thermo-electromagnetic model: cylindrical 3D domain and meridional section.

## 3.1 The electromagnetic model in Eulerian coordinates

Let us assume that the motion is given and, consequently, $\Omega(t)$ is known during the time interval $[0,T]$. $\Omega(t)$ will be composed by the main conducting elements present in the electric upsetting process: the cylindrical bar, the gripper jaws and the anvil (see Fig. 3). Let us assume that $\Omega(t)$ is cylindrically symmetric, namely

$$\Omega(t) := \left\{(r,\theta,z):\ \theta\in[0,2\pi),\ (r,z)\in\widehat{\Omega}(t)\right\},$$

for some bounded subset $\widehat{\Omega}(t) \subset \mathbb{R}^2$. Let $\mathbf{n}_{\mathsf{x}} = \widehat{n}_r \mathbf{e}_r + \widehat{n}_z \mathbf{e}_z$ be the outward unit normal vector to $\partial\Omega(t)$. We restrict our attention to a simply connected set $\Omega(t)$ that intersects the axis $r = 0$ in a set of positive one-dimensional measure, so that $\widehat{\Omega}(t)$ is also simply connected.

Let us further assume that the physical properties are independent on $\theta$ and that the current sources are such that the current density in the conducting part is of the form

$$\mathbf{J} := \widehat{J}_r \mathbf{e}_r + \widehat{J}_z \mathbf{e}_z \quad \text{in } \Omega(t).$$

Under these assumptions it is possible to define an eddy current model restricted to the conducting domain $\Omega(t)$ by using suitable boundary conditions. Firstly, we notice that if the current source is alternating and the materials have a linear magnetic behaviour, a time-harmonic approximation is often a suitable solution to avoid to work with time-scales which can be very different for thermal and electromagnetic phenomena; see, for instance, [18, 19]. However, even considering materials with a non-linear magnetic behaviour, this time-harmonic approximation is sometimes used to avoid long transient simulations which would be needed to reach a steady state in a genuine transient electromagnetic model; we refer the reader to [20] for further details in the complex representations of the ferromagnetic behaviour in non-linear time-harmonic problems. In this work, we adopt a time-harmonic approximation and consider the following eddy current model restricted to $\Omega(t)$:

$$i\omega\mathbf{B} + \mathbf{curl}\,\mathbf{E} = \mathbf{0} \quad \text{in } \Omega(t), \tag{23a}$$

$$\mathbf{curl}\,\mathbf{H} = \mathbf{J} \quad \text{in } \Omega(t), \tag{23b}$$

$$\operatorname{div}\mathbf{B} = 0 \quad \text{in } \Omega(t), \tag{23c}$$

$$\mathbf{B} = \breve{\mu}(|\mathbf{H}|, \Theta)\mathbf{H}, \tag{23d}$$

$$\mathbf{J} = \breve{\sigma}(\Theta)\mathbf{E}, \tag{23e}$$

where $\mathbf{B}$, $\mathbf{H}$, $\mathbf{E}$ are the complex amplitudes associated to the magnetic induction, the magnetic field and the electric field, respectively; $\omega$ is the angular frequency, i.e., $\omega = 2\pi f$ with $f$ the electric current frequency; $\breve{\mu}$ is the magnetic permeability which can be dependent on the temperature $\Theta$ and the modulus of the magnetic field $|\mathbf{H}|$; $\breve{\sigma}$ is the electrical conductivity, a function of temperature. We refer the reader to the book by Bossavit [15] for a discussion of the range of parameters for which the eddy current (or quasi-static) approximation is valid. This assumption is usually reasonable for low frequencies, as it happens in the electric upsetting process.

In the case of a fixed domain, the complex amplitudes do not depend on time. However, in the present case, since $\Omega$ changes with time, it is needed to consider the time-dependence in all the complex fields.

Next, we follow similar arguments as those developed in [16] to describe the problem in a bounded conducting domain under axisymmetric assumptions. Thus, by assuming that none of the components of the fields depends on $\theta$, we can look for a solution of the previous equations satisfying

$$\mathbf{E}(\mathsf{x}, t) := \widehat{E}_r(\widehat{\mathsf{x}}, t)\mathbf{e}_r + \widehat{E}_z(\widehat{\mathsf{x}}, t)\mathbf{e}_z \quad \text{in } \Omega(t),$$

$$\mathbf{H}(\mathrm{x}, t) := \widehat{H}_\theta(\widehat{\mathrm{x}}, t)\mathbf{e}_\theta \quad \text{in } \Omega(t).$$

Consequently, the following boundary condition can be imposed on the whole boundary of the conducting domain:

$$\check{\mu}\mathbf{H} \cdot \mathbf{n}_\mathrm{x} = 0 \quad \text{on } \partial\Omega(t). \tag{24}$$

This property will allow us to set boundary conditions on $\partial\Omega(t)$ that impose currents and/or potential drops on electric ports. Indeed, let us assume that the boundary of $\Omega(t)$ splits as follows:

$$\partial\Omega(t) := \Gamma_{\mathrm{N}}(t) \cup \Gamma_{\mathrm{J}}(t) \cup \Gamma_{\mathrm{E}}(t),$$

where $\Gamma_{\mathrm{J}}(t)$ and $\Gamma_{\mathrm{E}}(t)$ are the parts of the boundary connected to an electric source with known potential drops or currents, while $\Gamma_{\mathrm{N}}(t)$ is the isolated part, i.e., there is no current flux through this boundary. In eddy current models with electric ports it is usual to assume that currents enter and exit the domain perpendicularly and, consequently, we will assume

$$\mathbf{E} \times \mathbf{n}_\mathrm{x} = \mathbf{0} \quad \text{on } \Gamma_{\mathrm{J}}(t) \cup \Gamma_{\mathrm{E}}(t), \tag{25}$$

while the isolation condition means

$$\mathbf{J} \cdot \mathbf{n}_\mathrm{x} = \mathbf{curl}\,\mathbf{H} \cdot \mathbf{n}_\mathrm{x} = 0 \quad \text{on } \Gamma_{\mathrm{N}}(t). \tag{26}$$

From condition (24) we can deduce that there exists a sufficiently smooth function $\mathrm{U}(t)$ defined in $\Omega(t)$ up to a constant, such that $\mathrm{U}|_{\partial\Omega(t)}$ is a surface potential of the tangential component of $\mathbf{E}$, namely, $\mathbf{E} \times \mathbf{n}_\mathrm{x} = -\,\mathrm{grad}\,\mathrm{U} \times \mathbf{n}_\mathrm{x}$ on $\partial\Omega(t)$. On the other hand, (25) implies that $\mathrm{U}$ must be constant on each connected component of $\Gamma_{\mathrm{J}}(t) \cup \Gamma_{\mathrm{E}}(t)$ to be called a port. We assume that the whole $\Gamma_{\mathrm{E}}(t)$ is a port and denote the ports of $\Gamma_{\mathrm{J}}(t)$ as $\Gamma_{\mathrm{J}}^k(t)$, $k$ being its number. We can assume that $\mathrm{U} = 0$ on $\Gamma_{\mathrm{E}}(t)$ and then the complex number $\mathrm{U}_k := \mathrm{U}|_{\Gamma_{\mathrm{J}}^k(t)} - \mathrm{U}|_{\Gamma_{\mathrm{E}}(t)}$ is the potential drop between $\Gamma_{\mathrm{J}}^k(t)$ and $\Gamma_{\mathrm{E}}(t)$; consequently, $\mathrm{U}_k := \mathrm{U}|_{\Gamma_{\mathrm{J}}^k(t)}$.

According to the previous discussion, for each surface $\Gamma_{\mathrm{J}}^k(t)$, we will assume that one of the following is known:

- potential drop (or voltage) $\mathrm{U}_k := \mathrm{U}|_{\Gamma_{\mathrm{J}}^k(t)}$,
- current intensity through $\Gamma_{\mathrm{J}}^k(t)$, i.e.,

$$\int_{\Gamma_{\mathrm{J}}^k(t)} \mathbf{J} \cdot \mathbf{n}_\mathrm{x}\, \mathrm{d}A_\mathrm{x} = \int_{\Gamma_{\mathrm{J}}^k(t)} \mathbf{curl}\,\mathbf{H} \cdot \mathbf{n}_\mathrm{x}\, \mathrm{d}A_\mathrm{x} = I_k.$$

To obtain a weak formulation of this problem, let us multiply the Faraday's equation (23a) defined in $\Omega(t)$ by a smooth test vector field $\mathbf{G}$ such that $\mathbf{curl}\,\mathbf{G} \cdot \mathbf{n}_\mathrm{x} = 0$ on $\Gamma_{\mathrm{N}}(t)$. From Ampère's and Ohm's laws we obtain

$$\int_{\Omega(t)} i\omega\breve{\mu}\mathbf{H}\cdot\bar{\mathbf{G}}\,dV_x + \int_{\Omega(t)} \mathbf{E}\cdot\mathbf{curl}\,\bar{\mathbf{G}}\,dV_x$$
$$= \int_{\Omega(t)} i\omega\breve{\mu}\mathbf{H}\cdot\bar{\mathbf{G}}\,dV_x + \int_{\Omega(t)} \frac{1}{\breve{\sigma}}\mathbf{curl}\,\mathbf{H}\cdot\mathbf{curl}\,\bar{\mathbf{G}}\,dV_x = -\int_{\Omega(t)} \mathrm{grad}\,\mathrm{U}\cdot\mathbf{curl}\,\bar{\mathbf{G}}\,dV_x,$$

and using a Green's formula

$$\int_{\Omega(t)} i\omega\breve{\mu}\mathbf{H}\cdot\bar{\mathbf{G}}\,dV_x + \int_{\Omega(t)} \frac{1}{\breve{\sigma}}\mathbf{curl}\,\mathbf{H}\cdot\mathbf{curl}\,\bar{\mathbf{G}}\,dV_x = -\int_{\Omega(t)} \mathrm{grad}\,\mathrm{U}\cdot\mathbf{curl}\,\bar{\mathbf{G}}\,dV_x$$
$$= -\int_{\partial\Omega(t)} \mathrm{U}\,\mathbf{curl}\,\bar{\mathbf{G}}\cdot\mathbf{n}_x\,dA_x = -\int_{\Gamma_J(t)} \mathrm{U}\,\mathbf{curl}\,\bar{\mathbf{G}}\cdot\mathbf{n}_x\,dA_x, \quad (27)$$

where in the last equality we have used that $\mathrm{U}=0$ on $\Gamma_E(t)$.

We will distinguish between ports in $\Gamma_J(t)$ where we know the currents and those where the voltages are given; more precisely, the set of $N$ indices corresponding to the connected components of $\Gamma_J(t)$ is divided into two disjoint subsets: $\{1,\dots,N\} = N_I \cup N_V$, where

- For $k \in N_V$, voltage $V_k(t) \in \mathbb{C}$ is given.
- For $k \in N_I$, current $I_k(t) \in \mathbb{C}$ is given.

Fig. 4 shows the different parts of the boundary for a typical configuration of an upsetting machine. The electric ports connected to the source are $\Gamma_J^1$ on the lateral surface of the gripper and $\Gamma_E$ at the bottom of the anvil.

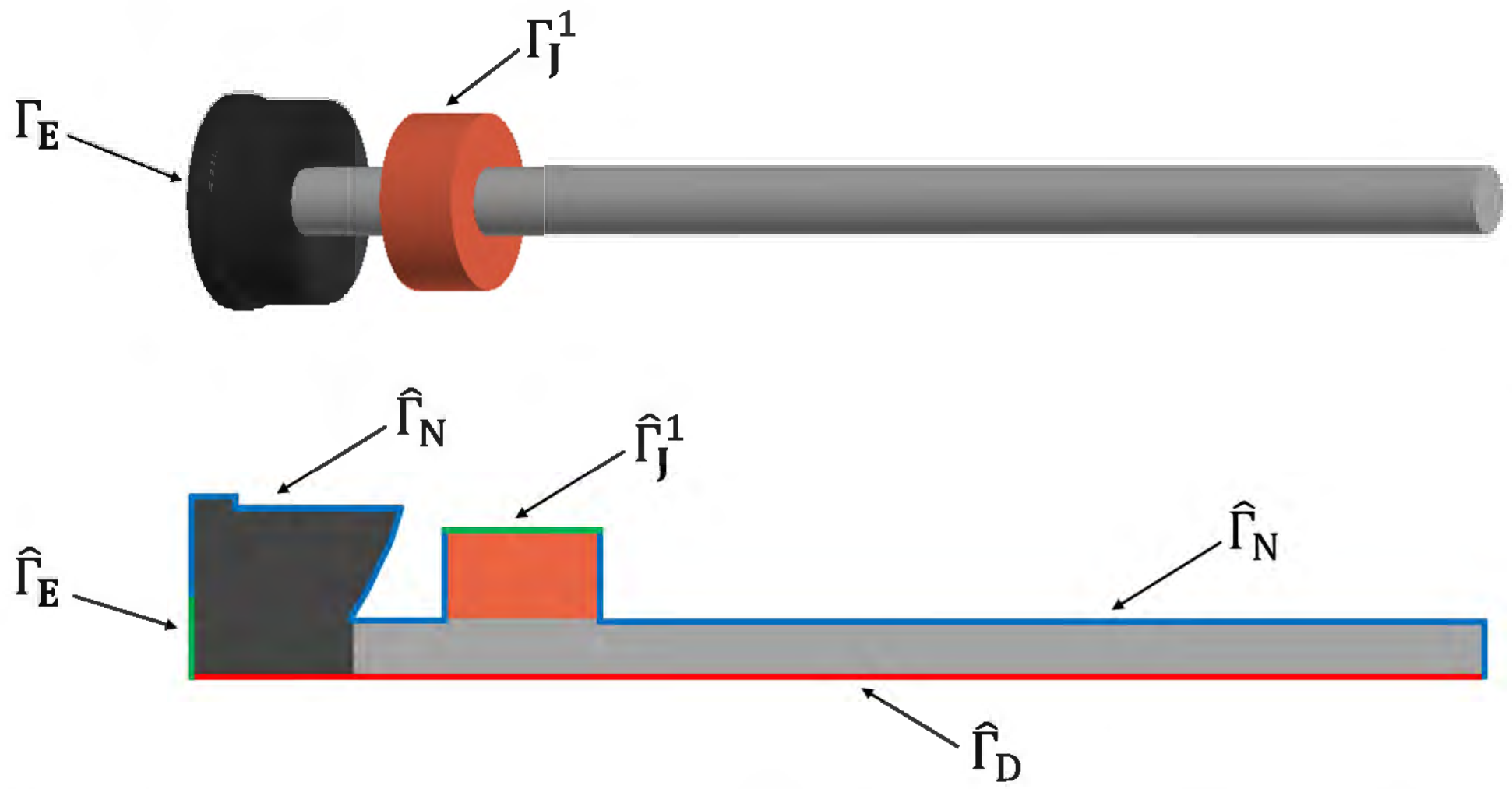

**Fig. 4**: Sketch of the conducting parts in an upsetting machine: 3D domain (above) and its meridional section (below).

Thus, the weak formulation of the electromagnetic problem can be written as follows:

*Given the displacement vector field* $\mathbf{u}(\mathsf{p},t)$*, the temperature* $\Theta(\mathsf{x},t)$*, and complex functions of time* $I_k(t)$ *for* $k\in N_{\rm I}$ *and* $\mathrm{V}_k(t)$ *for* $k\in N_{\rm V}$*, find complex fields* $\mathbf{H}(\mathsf{x},t)$*,* $\lambda(\mathsf{x},t)$ *defined on* $\Gamma_{\rm N}(t)$*, and* $\mathrm{V}_k(t)$ *for* $k\in N_{\rm I}$ *such that,*

$$\int_{\Omega(t)} i\omega\check{\mu}(|\mathbf{H}|,\Theta)\mathbf{H}\cdot\overline{\mathbf{G}}\,\mathrm{d}V_{\mathsf{x}}+\int_{\Omega(t)}\frac{1}{\check{\sigma}(\Theta)}\mathbf{curl}\mathbf{H}\cdot\mathbf{curl}\overline{\mathbf{G}}\,\mathrm{d}V_{\mathsf{x}}+\int_{\Gamma_{\rm N}(t)}\lambda\mathbf{curl}\overline{\mathbf{G}}\cdot\mathbf{n}_{\mathsf{x}}\,\mathrm{d}A_{\mathsf{x}}$$

$$+\sum_{k\in N_{\rm I}}\int_{\Gamma_{\rm J}^k(t)}\mathrm{V}_k\mathbf{curl}\overline{\mathbf{G}}\cdot\mathbf{n}_{\mathsf{x}}\,\mathrm{d}A_{\mathsf{x}}=-\sum_{k\in N_{\rm V}}\int_{\Gamma_{\rm J}^k(t)}\mathrm{V}_k\mathbf{curl}\overline{\mathbf{G}}\cdot\mathbf{n}_{\mathsf{x}}\,\mathrm{d}A_{\mathsf{x}}$$

$$\forall\mathbf{G}\text{ with }\mathbf{curl}\,\mathbf{G}\cdot\mathbf{n}_{\mathsf{x}}=0\text{ on }\Gamma_{\rm N}(t),$$

$$\sum_{k\in N_{\rm I}}\int_{\Gamma_{\rm J}^k(t)}\overline{\mathrm{W}}_k\mathbf{curl}\mathbf{H}\cdot\mathbf{n}_{\mathsf{x}}\,\mathrm{d}A_{\mathsf{x}}=\sum_{k\in N_{\rm I}}\overline{\mathrm{W}}_k\mathrm{I}_k\quad\forall\mathrm{W}_k,\ k\in N_{\rm I},$$

$$\int_{\Gamma_{\rm N}(t)}\overline{\eta}\mathbf{curl}\mathbf{H}\cdot\mathbf{n}_{\mathsf{x}}\,\mathrm{d}A_{\mathsf{x}}=0\quad\forall\eta.$$

Function $\lambda$ is defined on $\Gamma_{\rm N}(t)$ and is the Lagrange multiplier associated to the constraint of null current flux through the isolated boundaries, while $\mathrm{V}_k(t)$, $k\in N_{\rm I}$, are the Lagrange multipliers associated to the constraints imposing the currents in some electric ports.

Notice that, in general, for 3D domains it is not possible to restrict the model to conducting domains. In that case, a useful approach consists in using a magnetic field in the conducting parts and a scalar potential in the dielectric ones; see, for instance, [12] for a reference work by using these unknowns, applied later to the case of electric ports as in the present work, both in harmonic and transient regime (see, respectively [21, 22]).

Next, we will rewrite this weak formulation in cylindrical coordinates. Let us denote by $\partial\widehat{\Omega}(t)$ the boundary of $\widehat{\Omega}(t)$ which can be decomposed as $\partial\widehat{\Omega}(t):=\widehat{\Gamma}_{\rm D}(t)\cup\widehat{\Gamma}_{\rm J}(t)\cup\widehat{\Gamma}_{\rm N}(t)\cup\widehat{\Gamma}_{\rm E}(t)$ where we recall that $\widehat{\Gamma}_{\rm D}(t)$ is the part corresponding to the axis of symmetry $r=0$.

We consider a vector test function $\mathbf{G}(\mathsf{x},t)=\widehat{G}_\theta(\widehat{\mathsf{x}},t)\mathbf{e}_\theta$. Notice that

$$\mathbf{curl}\mathbf{G}\cdot\mathbf{n}_{\mathsf{x}}=\frac{1}{r}\widehat{\mathbf{grad}}(r\widehat{G}_\theta)\cdot\widehat{\boldsymbol{\tau}}_{\mathsf{x}}=\frac{1}{r}\frac{\partial(r\widehat{G}_\theta)}{\partial\widehat{\boldsymbol{\tau}}_{\mathsf{x}}}.\qquad(28)$$

Let us assume that $\widehat{\Gamma}_{\mathrm{J}}(t)$ has $N$ connected components such that $\widehat{\Gamma}_{\mathrm{J}}(t) := \bigcup_{k=1}^{N} \widehat{\Gamma}_{\mathrm{J}}^{k}(t)$. For $k = 1, \ldots, N$, we have,

$$\int_{\Gamma_{\mathrm{J}}^{k}(t)} \mathbf{curl}\,\mathbf{H}\cdot\mathbf{n}_{\mathrm{x}}\,\mathrm{d}A_{\mathrm{x}} = I_k = 2\pi\int_{\widehat{\Gamma}_{\mathrm{J}}^{k}(t)} \frac{1}{r}\frac{\partial(r\widehat{H}_\theta)}{\partial\widehat{\boldsymbol{\tau}}_{\mathrm{x}}} r\mathrm{d}\widehat{l}_{\mathrm{x}}$$

$$= 2\pi\int_{\widehat{\Gamma}_{\mathrm{J}}^{k}(t)} \frac{\partial(r\widehat{H}_\theta)}{\partial\widehat{\boldsymbol{\tau}}_{\mathrm{x}}}\mathrm{d}\widehat{l}_{\mathrm{x}}.$$

It is useful to make the change of variable $\tilde{H}_\theta := r\widehat{H}_\theta$ and $\tilde{G}_\theta := r\widehat{G}_\theta$. This leads to the following weak formulation written in the meridional section $\widehat{\Omega}(t)$:

*Given the displacement vector field $\widehat{\mathbf{u}}(\widehat{\mathsf{p}}, t)$, the temperature $\widehat{\Theta}(\widehat{\mathsf{x}}, t)$, and the complex-valued functions $I_k(t)$ for $k \in N_{\mathrm{I}}$ and $\mathrm{V}_k(t)$ for $k \in N_{\mathrm{V}}$, find complex fields $\tilde{H}_\theta(\widehat{\mathsf{x}}, t)$ being $\tilde{H}_\theta = 0$ on $\widehat{\Gamma}_{\mathrm{D}}(t)$, $\lambda(t)$ and $\mathrm{V}_k(t)$ for $k \in N_{\mathrm{I}}$, such that*

$$\int_{\widehat{\Omega}(t)} \frac{i\omega\check{\mu}(|\tilde{H}_\theta/r|, \widehat{\Theta})}{r}\tilde{H}_\theta\overline{\tilde{G}_\theta}\,\mathrm{d}r\mathrm{d}z + \int_{\widehat{\Omega}(t)} \frac{1}{\check{\sigma}(\widehat{\Theta})r}\left(\frac{\partial\tilde{H}_\theta}{\partial z}\frac{\partial\overline{\tilde{G}_\theta}}{\partial z} + \frac{\partial\tilde{H}_\theta}{\partial r}\frac{\partial\overline{\tilde{G}_\theta}}{\partial r}\right)\mathrm{d}r\mathrm{d}z$$

$$+ \int_{\widehat{\Gamma}_{\mathrm{N}}(t)} \lambda\frac{\partial\overline{\tilde{G}_\theta}}{\partial\widehat{\boldsymbol{\tau}}_{\mathrm{x}}}\mathrm{d}\widehat{l}_{\mathrm{x}} + \sum_{k\in N_{\mathrm{I}}}\int_{\widehat{\Gamma}_{\mathrm{J}}^{k}(t)} \mathrm{V}_k\frac{\partial\overline{\tilde{G}_\theta}}{\partial\widehat{\boldsymbol{\tau}}_{\mathrm{x}}}\mathrm{d}\widehat{l}_{\mathrm{x}}$$

$$= -\sum_{k\in N_{\mathrm{V}}}\int_{\widehat{\Gamma}_{\mathrm{J}}^{k}(t)} \mathrm{V}_k\frac{\partial\overline{\tilde{G}_\theta}}{\partial\widehat{\boldsymbol{\tau}}_{\mathrm{x}}}\mathrm{d}\widehat{l}_{\mathrm{x}} \quad \forall\tilde{G}_\theta \text{ with } \tilde{G}_\theta = 0 \text{ en } \widehat{\Gamma}_{\mathrm{D}}(t),$$

$$\sum_{k\in N_{\mathrm{I}}}\int_{\widehat{\Gamma}_{\mathrm{J}}^{k}(t)} \overline{\mathrm{W}}_k\frac{\partial\tilde{H}_\theta}{\partial\widehat{\boldsymbol{\tau}}_{\mathrm{x}}}\mathrm{d}\widehat{l}_{\mathrm{x}} = \sum_{k\in N_{\mathrm{I}}}\overline{\mathrm{W}}_k\frac{\mathrm{I}_k}{2\pi} \quad \forall\mathrm{W}_k,\ k\in N_{\mathrm{I}},$$

$$\int_{\widehat{\Gamma}_{\mathrm{N}}(t)} \overline{\eta}\,\frac{\partial\tilde{H}_\theta}{\partial\widehat{\boldsymbol{\tau}}_{\mathrm{x}}}\mathrm{d}\widehat{l}_{\mathrm{x}} = 0 \quad \forall\eta.$$

## 3.2 A Lagrangian approach to the eddy current model

From the weak formulation in Eulerian coordinates presented above, we will obtain a weak formulation in Lagrangian coordinates. First, we will perform the computations in the 3D case to later obtain the axisymmetric formulation.

First, notice that the **curl** operator for an Eulerian vector field, $\boldsymbol{\Psi}(\mathsf{x}, t)$, can be written in Lagrangian coordinates as (see [23] for further details)

$$\mathbf{curl}\boldsymbol{\Psi}(\mathsf{x}, t) = \frac{1}{\det\mathbf{F}(\mathsf{p}, t)}\mathbf{F}(\mathsf{p}, t)\,\mathbf{Curl}\left(\mathbf{F}^t(\mathsf{p}, t)\,\boldsymbol{\Psi}_m(\mathsf{p}, t)\right)\Big|_{\mathsf{p}=\mathsf{P}(\mathsf{x},t)}, \tag{29}$$

being $\boldsymbol{\Psi}_m$ the material description of the vector field $\boldsymbol{\Psi}$.

This property leads us to introduce a new field $\boldsymbol{\mathcal{H}}(\mathsf{p},t)$ which is related to the material description of the unknown $\mathbf{H}$, $\mathbf{H}_m(\mathsf{p},t)$, as follows:

$$\boldsymbol{\mathcal{H}}(\mathsf{p},t) = \mathbf{F}^t(\mathsf{p},t)\,\mathbf{H}_m(\mathsf{p},t). \tag{30}$$

In a similar way, a new function $\boldsymbol{\mathcal{G}}(\mathsf{p},t)$ associated with the test function $\mathbf{G}$ will be introduced, satisfying a similar relation to (30).

By using (29) and (30) and identities (8)-(9), we obtain a 3D weak formulation written in Lagrangian coordinates:

*Given the displacement vector field* $\mathbf{u}(\mathsf{p},t)$*, the temperature* $\Theta_m(\mathsf{p},t)$*, and the complex-valued functions* $I_k(t)$ *for* $k \in N_{\mathrm{I}}$ *and* $\mathrm{V}_k(t)$ *for* $k \in N_{\mathrm{V}}$*, find complex fields* $\boldsymbol{\mathcal{H}}(\mathsf{p},t)$*,* $\lambda(t)$ *and* $\mathrm{V}_k(t)$ *for* $k \in N_{\mathrm{I}}$ *such that*

$$\begin{aligned}
\int_\Omega i\omega\check{\mu}(|\mathbf{F}^{-t}\boldsymbol{\mathcal{H}}|,\Theta_m)\mathbf{F}^{-t}\boldsymbol{\mathcal{H}}\cdot\mathbf{F}^{-t}\overline{\boldsymbol{\mathcal{G}}}\det\mathbf{F}\,\mathrm{d}V_{\mathsf{p}} + \int_\Omega \frac{1}{\check{\sigma}(\Theta_m)\det\mathbf{F}}\mathbf{F}\mathbf{Curl}\boldsymbol{\mathcal{H}}\cdot\mathbf{F}\mathbf{Curl}\overline{\boldsymbol{\mathcal{G}}}\,\mathrm{d}V_{\mathsf{p}}& \\
+\int_{\Gamma_{\mathrm{N}}}\lambda\mathbf{Curl}\overline{\boldsymbol{\mathcal{G}}}\cdot\mathbf{n}_{\mathsf{p}}\,\mathrm{d}A_{\mathsf{p}} + \sum_{k\in N_{\mathrm{I}}}\int_{\Gamma_{\mathrm{J}}^k}\mathrm{V}_k\mathbf{Curl}\overline{\boldsymbol{\mathcal{G}}}\cdot\mathbf{n}_{\mathsf{p}}\,\mathrm{d}A_{\mathsf{p}}& \\
= -\sum_{k\in N_{\mathrm{V}}}\int_{\Gamma_{\mathrm{J}}^k}\mathrm{V}_k\mathbf{Curl}\overline{\boldsymbol{\mathcal{G}}}\cdot\mathbf{n}_{\mathsf{p}}\,\mathrm{d}A_{\mathsf{p}}\ \forall\boldsymbol{\mathcal{G}}\text{ with }\mathbf{Curl}\boldsymbol{\mathcal{G}}\cdot\mathbf{n}_{\mathsf{p}} = 0\text{ on }\Gamma_{\mathrm{N}},& \\
\sum_{k\in N_{\mathrm{I}}}\int_{\Gamma_{\mathrm{J}}^k}\overline{\mathrm{W}}_k\mathbf{Curl}\boldsymbol{\mathcal{H}}\cdot\mathbf{n}_{\mathsf{p}}\,\mathrm{d}A_{\mathsf{p}} = \sum_{k\in N_{\mathrm{I}}}\overline{\mathrm{W}}_k\mathrm{I}_k\quad\forall\mathrm{W}_k,\ k\in N_{\mathrm{I}},& \\
\int_{\Gamma_{\mathrm{N}}}\overline{\eta}\mathbf{Curl}\boldsymbol{\mathcal{H}}\cdot\mathbf{n}_{\mathsf{p}}\,\mathrm{d}A_{\mathsf{p}} = 0\quad\forall\eta.&
\end{aligned}$$

By taking into account the cylindrical symmetry, as we did in the Eulerian case, we consider

$$\begin{aligned}
\boldsymbol{\mathcal{H}}(\mathsf{p},t) &= \widehat{\mathcal{H}}_\theta(\widehat{\mathsf{p}},t)\mathbf{e}_\theta\text{ in }\Omega,\\
\boldsymbol{\mathcal{G}}(\mathsf{p},t) &= \widehat{\mathcal{G}}_\theta(\widehat{\mathsf{p}},t)\mathbf{e}_\theta\text{ in }\Omega.
\end{aligned}$$

Thus, by introducing the change of variable $\tilde{\mathcal{H}}_\theta = r_m\widehat{\mathcal{H}}_\theta$ and $\tilde{\mathcal{G}}_\theta = r_m\widehat{\mathcal{G}}_\theta$, we get the following relations which will be used to obtain the weak formulation in the axisymmetric case:

$$\begin{aligned}
\mathbf{F}^{-t}(\mathsf{p},t)\,\boldsymbol{\mathcal{H}}(\mathsf{p},t)\cdot\mathbf{F}^{-t}(\mathsf{p},t)\,\overline{\boldsymbol{\mathcal{G}}}(\mathsf{p},t) &= \frac{1}{(r_m+\widehat{u}_r(\widehat{\mathsf{p}},t))^2}\tilde{\mathcal{H}}_\theta(\widehat{\mathsf{p}},t)\overline{\tilde{\mathcal{G}}_\theta}(\widehat{\mathsf{p}},t),\\
\mathbf{F}(\mathsf{p},t)\,\mathbf{Curl}\boldsymbol{\mathcal{H}}(\mathsf{p},t)\cdot\mathbf{F}(\mathsf{p},t)\,\mathbf{Curl}\overline{\boldsymbol{\mathcal{G}}}(\mathsf{p},t) &= \\
\frac{1}{r_m^2}\widehat{\widehat{\mathbf{N}}}(\widehat{\mathsf{p}},t)\,\widehat{\mathbf{Grad}\tilde{\mathcal{H}}_\theta}(\widehat{\mathsf{p}},t)&\cdot\widehat{\widehat{\mathbf{N}}}(\widehat{\mathsf{p}},t)\,\widehat{\mathbf{Grad}\overline{\tilde{\mathcal{G}}_\theta}}(\widehat{\mathsf{p}},t),
\end{aligned}$$

being $\widehat{\widehat{\mathbf{N}}}$ the tensor field

$$\widehat{\widehat{\mathbf{N}}}(\widehat{\mathsf{p}},t) = \begin{pmatrix} \dfrac{\partial \widehat{u}_r(\widehat{\mathsf{p}},t)}{\partial z_m} & -1 - \dfrac{\partial \widehat{u}_r(\widehat{\mathsf{p}},t)}{\partial r_m} \\ 1 + \dfrac{\partial \widehat{u}_z(\widehat{\mathsf{p}},t)}{\partial z_m} & -\dfrac{\partial \widehat{u}_z(\widehat{\mathsf{p}},t)}{\partial r_m} \end{pmatrix}. \tag{31}$$

Finally, by using (21), the Lagrangian weak formulation defined in the meridional section $\widehat{\Omega}$ reads as follows:

*Given the displacement vector field $\widehat{\mathbf{u}}(\widehat{\mathsf{p}},t)$, the temperature $\widehat{\Theta}_m(\widehat{\mathsf{p}},t)$, and complex-valued functions $I_k(t)$ for $k \in N_{\rm I}$ and $\mathrm{V}_k(t)$ for $k \in N_{\rm V}$, find complex fields $\tilde{\mathcal{H}}_\theta(\widehat{\mathsf{p}},t)$, with $\tilde{\mathcal{H}}_\theta = 0$ on $\widehat{\Gamma}_{\rm D}$, $\lambda(t)$ and $\mathrm{V}_k(t)$ for $k \in N_{\rm I}$ such that*

$$\int_{\widehat{\Omega}} \frac{i\omega\breve{\mu}(|\tilde{\mathcal{H}}_\theta/(r_m+\widehat{u}_r)|,\widehat{\Theta}_m)}{r_m+\widehat{u}_r}\tilde{\mathcal{H}}_\theta\overline{\tilde{\mathcal{G}}_\theta}\,\det\widehat{\widehat{\mathbf{F}}}\,\mathrm{d}r_m\mathrm{d}z_m$$

$$+\int_{\widehat{\Omega}} \frac{\widehat{\widehat{\mathbf{N}}}\,\widehat{\mathbf{Grad}\tilde{\mathcal{H}}_\theta}\cdot\widehat{\widehat{\mathbf{N}}}\,\widehat{\mathbf{Grad}\overline{\tilde{\mathcal{G}}_\theta}}}{\breve{\sigma}(\widehat{\Theta}_m)\,(r_m+\widehat{u}_r)\,\det\widehat{\widehat{\mathbf{F}}}}\mathrm{d}r_m\mathrm{d}z_m + \int_{\widehat{\Gamma}_{\rm N}} \lambda\frac{\partial\overline{\tilde{\mathcal{G}}_\theta}}{\partial\widehat{\boldsymbol{\tau}}_{\mathsf{p}}}\mathrm{d}\widehat{l}_{\mathsf{p}} + \sum_{k\in N_{\rm I}}\int_{\widehat{\Gamma}_{\rm J}^k} \mathrm{V}_k\frac{\partial\overline{\tilde{\mathcal{G}}_\theta}}{\partial\widehat{\boldsymbol{\tau}}_{\mathsf{p}}}\mathrm{d}\widehat{l}_{\mathsf{p}}$$

$$= -\sum_{k\in N_{\rm V}}\int_{\widehat{\Gamma}_{\rm J}^k} \mathrm{V}_k\frac{\partial\overline{\tilde{\mathcal{G}}_\theta}}{\partial\widehat{\boldsymbol{\tau}}_{\mathsf{p}}}\mathrm{d}\widehat{l}_{\mathsf{p}} \quad \forall\tilde{\mathcal{G}}_\theta \text{ with } \tilde{\mathcal{G}}_\theta = 0 \text{ on } \widehat{\Gamma}_{\rm D},$$

$$\sum_{k\in N_{\rm I}}\int_{\widehat{\Gamma}_{\rm J}^k} \overline{\mathrm{W}}_k\frac{\partial\tilde{\mathcal{H}}_\theta}{\partial\widehat{\boldsymbol{\tau}}_{\mathsf{p}}}\mathrm{d}\widehat{l}_{\mathsf{p}} = \sum_{k\in N_{\rm I}}\overline{\mathrm{W}}_k\frac{\mathrm{I}_k}{2\pi} \quad \forall\mathrm{W}_k,\ k\in N_{\rm I},$$

$$\int_{\widehat{\Gamma}_{\rm N}} \overline{\eta}\,\frac{\partial\tilde{\mathcal{H}}_\theta}{\partial\widehat{\boldsymbol{\tau}}_{\mathsf{p}}}\mathrm{d}\widehat{l}_{\mathsf{p}} = 0 \quad \forall\eta.$$

# 4 The thermal model

In this section we introduce the thermal model which allows us to compute the temperature in a meridional section of the cylindrical domain $\Omega(t)$ due to the power dissipated by Joule effect in the electromagnetic model. We follow the same scheme as in the electromagnetic case: first, we introduce the problem in 3D and Eulerian coordinates to finally get the Lagrangian formulation in the axisymmetric setting.

## 4.1 The thermal model in Eulerian coordinates

Let us consider the same domain $\Omega(t)$ of the electromagnetic model. In this case, the boundary $\Gamma(t)$ is decomposed in a part $\Gamma_{\rm DT}(t)$ where the temperature is known, and a part $\Gamma_{\rm CR}(t)$ where convection-radiation boundary conditions are imposed. Thus, $\partial\Omega(t) = \Gamma_{\rm DT}(t) \cup \Gamma_{\rm CR}(t)$. The heat transfer problem in transient state and Eulerian coordinates consists in finding a temperature field $\Theta(\mathsf{x},t)$ satisfying

$$\rho\check{c}_p(\Theta)\left(\frac{\partial\Theta}{\partial t}+\mathbf{v}\cdot\mathbf{grad}\Theta\right)-\operatorname{div}(\check{k}(\Theta)\mathbf{grad}\Theta)=\frac{|\mathbf{J}|^2}{2\check{\sigma}(\Theta)}\quad\text{in }\Omega(t), \tag{32}$$

$$\Theta=\Theta^{\mathrm{DT}}\quad\text{on }\Gamma_{\mathrm{DT}}(t), \tag{33}$$

$$\check{k}(\Theta)\frac{\partial\Theta}{\partial\mathbf{n}_{\mathsf{x}}}=h(\Theta_C-\Theta)+\sigma_{SB}\epsilon(\Theta_R^4-\Theta^4)\quad\text{on }\Gamma_{\mathrm{CR}}(t), \tag{34}$$

$$\Theta(:,0)=\Theta^0\quad\text{in }\bar{\Omega}, \tag{35}$$

where $\check{c}_p$ and $\check{k}$ are, respectively, the specific heat and the thermal conductivity as functions of temperature; $\rho$ is the mass density; $h$ is the heat transfer coefficient; $\sigma_{SB}$ is the Stefan-Boltzmann constant; $\epsilon$ is the emissivity; $\Theta_C$ is the convection temperature, $\Theta_R$ is the radiation temperature and $\Theta^0$ is the initial temperature field. The source term is the Joule effect which is computed from the electromagnetic solution and couples the thermal and the electromagnetic models.

The weak formulation of this problem is standard:

*Given the displacement vector field $\mathbf{u}(\mathsf{p},t)$, the temperature field at initial time, $\Theta^0(\mathsf{p})$, the temperature $\Theta^{\mathrm{DT}}(\mathsf{x},t)$ on $\Gamma_{\mathrm{DT}}(t)$ and the modulus of the current density, $|\mathbf{J}(\mathsf{x},t)|$, find a scalar field $\Theta(\mathsf{x},t)$ such that $\Theta(\mathsf{x},t)=\Theta^{\mathrm{DT}}(\mathsf{x},t)$ on $\Gamma_{\mathrm{DT}}(t)$ and $\Theta(:,0)=\Theta^0$ in $\bar{\Omega}$ satisfying*

$$\int_{\Omega(t)}\rho\,\check{c}_p(\Theta)\left(\frac{\partial\Theta}{\partial t}+\mathbf{v}\cdot\mathbf{grad}\Theta\right)\psi\,\mathrm{d}V_{\mathsf{x}}+\int_{\Omega(t)}\check{k}(\Theta)\,\mathbf{grad}\,\Theta\cdot\mathbf{grad}\,\psi\,\mathrm{d}V_{\mathsf{x}}$$
$$-\int_{\Gamma_{\mathrm{CR}}(t)}\left[h(\Theta_C-\Theta)+\sigma_{SB}\epsilon(\Theta_R^4-\Theta^4)\right]\psi\,\mathrm{d}A_{\mathsf{x}}=\int_{\Omega(t)}\frac{|\mathbf{J}|^2}{2\check{\sigma}(\Theta)}\psi\,\mathrm{d}V_{\mathsf{x}}$$
$$\forall\psi\text{ null on }\Gamma_{\mathrm{DT}}(t). \tag{36}$$

Now, we assume that none of the fields depends on the azimuthal coordinate, and consider the boundaries $\widehat{\Gamma}_{\mathrm{D}}(t)$, $\widehat{\Gamma}_{\mathrm{DT}}(t)$ and $\widehat{\Gamma}_{\mathrm{CR}}(t)$, where $\widehat{\Gamma}_{\mathrm{D}}(t)$ denotes the symmetry axis as in the previous section; notice that the heat flux will be null through $\widehat{\Gamma}_{\mathrm{D}}(t)$.

Thus, by using the notation in Section 2 applied to the thermal problem, and introducing

$$\widehat{\mathbf{v}}(\widehat{\mathsf{x}},t)=\widehat{v}_r(\widehat{\mathsf{x}},t)\widehat{\mathbf{e}}_r+\widehat{v}_z(\widehat{\mathsf{x}},t)\widehat{\mathbf{e}}_z\text{ in }\widehat{\Omega}(t), \tag{37}$$

we obtain its axisymmetric Eulerian weak formulation:

*Given the displacement vector field $\widehat{\mathbf{u}}(\widehat{\mathsf{p}},t)$, the temperature field at initial time, $\widehat{\Theta}^0(\widehat{\mathsf{p}})$, the temperature $\widehat{\Theta}^{\mathrm{DT}}(\widehat{\mathsf{x}},t)$ on $\widehat{\Gamma}_{\mathrm{DT}}(t)$ and the modulus of the current density,*

$|\widehat{\mathbf{J}}(\widehat{\mathsf{x}},t)|$, *find a scalar field,* $\widehat{\Theta}(\widehat{\mathsf{x}},t)$ *with* $\widehat{\Theta}(\widehat{\mathsf{x}},t) = \widehat{\Theta}^{\mathrm{DT}}(\widehat{\mathsf{x}},t)$ *on* $\widehat{\Gamma}_{\mathrm{DT}}(t)$ *and* $\widehat{\Theta}(:,0) = \widehat{\Theta}^0$ *in* $\widehat{\bar{\Omega}}$ *such that*

$$\int_{\widehat{\Omega}(t)} \widehat{\rho}\,\check{c}_p(\widehat{\Theta}) \left( \frac{\partial \widehat{\Theta}}{\partial t} + \widehat{\mathbf{v}} \cdot \widehat{\mathbf{grad}\Theta} \right) \widehat{\psi} r \,\mathrm{d}r\mathrm{d}z + \int_{\widehat{\Omega}(t)} \check{k}(\widehat{\Theta})\, \widehat{\mathbf{grad}\Theta} \cdot \widehat{\mathbf{grad}\psi} r \,\mathrm{d}r\mathrm{d}z$$
$$- \int_{\widehat{\Gamma}_{\mathrm{CR}}(t)} \left[ h(\widehat{\Theta}_C - \widehat{\Theta}) + \sigma_{SB}\epsilon(\widehat{\Theta}_R^4 - \widehat{\Theta}^4) \right] \widehat{\psi} r \,\mathrm{d}\widehat{l}_{\mathsf{x}} = \int_{\widehat{\Omega}(t)} \frac{|\widehat{\mathbf{J}}|^2}{2\check{\sigma}(\widehat{\Theta})} \widehat{\psi} r \,\mathrm{d}r\mathrm{d}z$$
$$\forall \widehat{\psi} \text{ null on } \widehat{\Gamma}_{\mathrm{DT}}(t). \quad (38)$$

## 4.2 A Lagrangian approach to the thermal model

By using the relation:

$$\mathbf{grad}\phi(\mathsf{x},t) = \mathbf{F}^{-t}(\mathsf{p},t)\, \mathbf{Grad}\phi_m(\mathsf{p},t)\big|_{\mathsf{p}=\mathsf{P}(\mathsf{x},t)}, \quad (39)$$

and the identities (6)-(8), we obtain the standard weak formulation of (32)-(35) in 3D Lagrangian coordinates:

*Given the displacement vector field* $\mathbf{u}(\mathsf{p},t)$, *the initial field temperature,* $\Theta^0(\mathsf{p})$, *the temperature* $\Theta_m^{\mathrm{DT}}(\mathsf{p},t)$ *on* $\Gamma_{\mathrm{DT}}$ *and the modulus of the current density,* $|\mathbf{J}_m(\mathsf{p},t)|$, *find a scalar field* $\Theta_m(\mathsf{p},t)$ *with* $\Theta_m(\mathsf{p},t) = \Theta_m^{\mathrm{DT}}(\mathsf{p},t)$ *on* $\Gamma_{\mathrm{DT}}$ *and* $\Theta_m(:,0) = \Theta^0$ *in* $\bar{\Omega}$ *such that*

$$\int_\Omega \rho_0\, \check{c}_p(\Theta_m)\dot{\Theta}_m \psi_m \,\mathrm{d}V_{\mathsf{p}} + \int_\Omega \check{k}(\Theta_m)\, \mathbf{F}^{-t}\mathbf{Grad}\Theta_m \cdot \mathbf{F}^{-t}\mathbf{Grad}\psi_m \det\mathbf{F}\,\mathrm{d}V_{\mathsf{p}}$$
$$- \int_{\Gamma_{\mathrm{CR}}} \left[ h(\Theta_C - \Theta_m) + \sigma_{SB}\epsilon(\Theta_R^4 - \Theta_m^4) \right] \psi_m |\mathbf{F}^{-\mathbf{t}}\mathbf{n}_{\mathsf{p}}| \det\mathbf{F}\,\mathrm{d}A_{\mathsf{p}}$$
$$= \int_\Omega \frac{|\mathbf{J}_m|^2}{2\check{\sigma}(\Theta_m)} \psi_m \det\mathbf{F}\,\mathrm{d}V_{\mathsf{p}} \quad \forall \psi_m \text{ null on } \Gamma_{\mathrm{DT}}. \quad (40)$$

Notice that we have used the mass conservation principle which states

$$\rho_m(\mathsf{p},t)\det\mathbf{F}(\mathsf{p},t) = \rho_0(\mathsf{p}), \quad (41)$$

being $\rho_0(\mathsf{p})$ the mass density in the reference configuration.

In the case of cylindrical symmetry, we use the notation given in Section 2 to reach the following axisymmetric Lagrangian weak formulation for the thermal model:

*Given the displacement vector field* $\widehat{\mathbf{u}}(\widehat{\mathsf{p}},t)$, *the initial field temperature,* $\widehat{\Theta}^0(\widehat{\mathsf{p}})$, *the temperature* $\widehat{\Theta}_m^{\mathrm{DT}}(\widehat{\mathsf{p}},t)$ *on* $\widehat{\Gamma}_{\mathrm{DT}}$ *and the modulus of the current density,* $|\widehat{\mathbf{J}}_m(\widehat{\mathsf{p}},t)|$, *find a scalar field,* $\widehat{\Theta}_m(\widehat{\mathsf{p}},t)$ *with* $\widehat{\Theta}_m(\widehat{\mathsf{p}},t) = \widehat{\Theta}_m^{\mathrm{DT}}(\widehat{\mathsf{p}},t)$ *on* $\widehat{\Gamma}_{\mathrm{DT}}$ *and* $\widehat{\Theta}_m(:,0) = \widehat{\Theta}^0$ *in* $\widehat{\bar{\Omega}}$ *such that*

$$
\begin{aligned}
&\int_{\widehat{\Omega}} \widehat{\rho}_0 \check{c}_p(\widehat{\Theta}_m)\dot{\widehat{\Theta}}_m \widehat{\psi}_m r_m \,\mathrm{d}r_m \mathrm{d}z_m \\
&+ \int_{\widehat{\Omega}} \check{k}(\widehat{\Theta}_m)\, \widehat{\widehat{\mathbf{F}}}^{-t}\, \widehat{\mathbf{Grad}\Theta}_m \cdot \widehat{\widehat{\mathbf{F}}}^{-t}\, \widehat{\mathbf{Grad}\psi}_m \left(1 + \frac{\widehat{u}_r}{r_m}\right) \det\widehat{\widehat{\mathbf{F}}} r_m \,\mathrm{d}r_m \mathrm{d}z_m \\
&- \int_{\widehat{\Gamma}_{\mathrm{CR}}} \left[h(\widehat{\Theta}_C - \widehat{\Theta}_m) + \sigma_{SB}\epsilon(\widehat{\Theta}_R^4 - \widehat{\Theta}_m^4)\right] \widehat{\psi}_m |\widehat{\widehat{\mathbf{F}}}^{-t} \widehat{\mathbf{n}}_{\mathrm{p}}| \left(1 + \frac{\widehat{u}_r}{r_m}\right) \det\widehat{\widehat{\mathbf{F}}} r_m \,\mathrm{d}\widehat{l}_{\mathrm{p}} \\
&= \int_{\widehat{\Omega}} \frac{|\widehat{\mathbf{J}}_m|^2}{2\check{\sigma}(\widehat{\Theta}_m)} \widehat{\psi}_m \left(1 + \frac{\widehat{u}_r}{r_m}\right) \det\widehat{\widehat{\mathbf{F}}} r_m \,\mathrm{d}r_m \mathrm{d}z_m \\
&\forall \widehat{\psi}_m \text{ null on } \widehat{\Gamma}_{\mathrm{DT}}. \qquad (42)
\end{aligned}
$$

# 5 Numerical solution of the coupled problem

The electromagnetic and thermal axisymmetric problems in Lagrangian coordinates have been discretized by using a time discretization and a finite element spatial discretization, and then implemented in Python-FEniCS [24]. For the time discretization we have used an implicit Euler method. Concerning the spatial discretization, the temperature and the electromagnetic variable $\widetilde{H}_\theta$ have been approximated by piecewise linear elements, and the unknown $\lambda$ by piecewise constant elements defined on the edges of the boundary $\widehat{\Gamma}_{\mathrm{N}}$.

To solve the coupled thermo-electromagnetic problem, we have used a monolithic scheme, i.e. a fully-coupled scheme in which the governing equations are solved simultaneously. Thus, at each time step, a non-linear problem involving the electromagnetic and thermal unknowns is solved using a Newton-Raphson method.

The numerical procedure has been validated by solving a thermo-electromagnetic test stated in a conducting cylinder of length $L = 0.165\,\mathrm{m}$ and radius $R = 0.02875\,\mathrm{m}$ that is deformed by a given displacement field that tries to emulate a typical electric upsetting. Let us consider a displacement vector field with a zero axial component and a radial component, $\widehat{u}_r(r_m, z_m)$, given by

$$
\widehat{u}_r(r_m, z_m) = \begin{cases} 10^3 r_m\,(-188.2593 z_m + 6.1464)\, z_m^2 & \text{if } z_m \leq 0.02, \\ r_m \left[1.0793 \exp\left(-\left(\dfrac{z_m - 0.0293}{0.03104}\right)^2\right) + \right. \\ \quad - 18.4974 \exp\left(-\left(\dfrac{z_m + 0.03324}{0.01705}\right)^2\right) + \\ \quad \left. +1.0779 \exp\left(-\left(\dfrac{z_m - 0.4363}{1.263}\right)^2\right) - 1\right] & \text{if } z_m > 0.02. \end{cases}
$$

An alternating current flows through the cylinder from top to bottom. The lateral surface of the cylinder is isolated, that is, the current flux is null. The input current is provided at the top of the cylinder, while the ground $\Gamma_{\mathbf{E}}$ is at the bottom; the current enters and exits normal to the electric ports as we stated in the model presented previously. The amplitude of the current is equal to 35000 A and the electrical frequency is 500 Hz.

The electromagnetic and thermal physical properties are temperature-dependent functions (with temperature $\Theta$ in °C) that are provided to the solver. These functions are detailed below and have been obtained by fitting tables of data corresponding to a typical steel. Figure 5a and Figure 5b show each property vs. temperature.

More specifically, the expressions for $\breve{\sigma}$, $\breve{k}$ and $\breve{c}_p$ are the following:

$$\begin{aligned}
\breve{\sigma}(\Theta) &= \frac{1}{-4.3306\times 10^{-13}\Theta^2 + 1.0839\times 10^{-9}\Theta + 2.0170\times 10^{-7}},\\
\breve{k}(\Theta) &= -2.7834\times 10^{-11}\Theta^4 + 1.1045\times 10^{-7}\Theta^3 - 1.3658\times 10^{-4}\Theta^2 + \\
&\quad + 0.04639\Theta + 34.0140,\\
\breve{c}_p(\Theta) &= 660.9\exp\left(-\left(\frac{\Theta - 723.3}{23.93}\right)^2\right) + 288.9\exp\left(-\left(\frac{\Theta - 697.6}{133.5}\right)^2\right) + \\
&\quad + 657.1\exp\left(-\left(\frac{\Theta - 908.1}{1497.0}\right)^2\right).
\end{aligned}$$

On the other hand, the magnetic permeability is a function of the temperature and also of the magnetic field, and is given by a modified version of the so-called Fröhlich–Kennelly model, described in [25]. More precisely, $\breve{\mu}$ in terms of the modulus of the magnetic field $|\mathbf{H}|$ and the temperature $\Theta$ is given by

$$\breve{\mu}(|\mathbf{H}|, \Theta) = \mu_0 + \frac{f(\Theta)}{a + b|\mathbf{H}|},$$

where $a = 2532.35\,\mathrm{T}^{-1}\mathrm{Am}^{-1}$ and $b = 0.49\,\mathrm{T}^{-1}$; $f(\Theta)$ is a temperature-dependent factor given by

$$f(\Theta) = \begin{cases} \left(\dfrac{(\Theta^C + 273.15)^2 - (\Theta + 273.15)^2}{(\Theta^C + 273.15)^2 - (\Theta^0 + 273.15)^2}\right)^{\frac{1}{4}} & \text{if } \Theta < \Theta^C, \\ 0 & \text{if } \Theta \geq \Theta^C, \end{cases}$$

where $\Theta^C = 748.69$ °C is the Curie temperature of the material and $\Theta^0 = 23.5$ °C is the room temperature. The relationship between $H_\theta$ and $B_\theta$ is shown in Figure 5a-left for different values of temperature.

The mass density is constant and equal to $\rho_m = 7799\ \mathrm{kg\,m^{-3}}$. Finally, the initial temperature of the cylinder is equal to 20 °C and there is no heat flux through the boundaries.

To validate the code in Lagrangian coordinates, we also solved the problem in Eulerian coordinates with a similar discretization and compared the results. The coupled

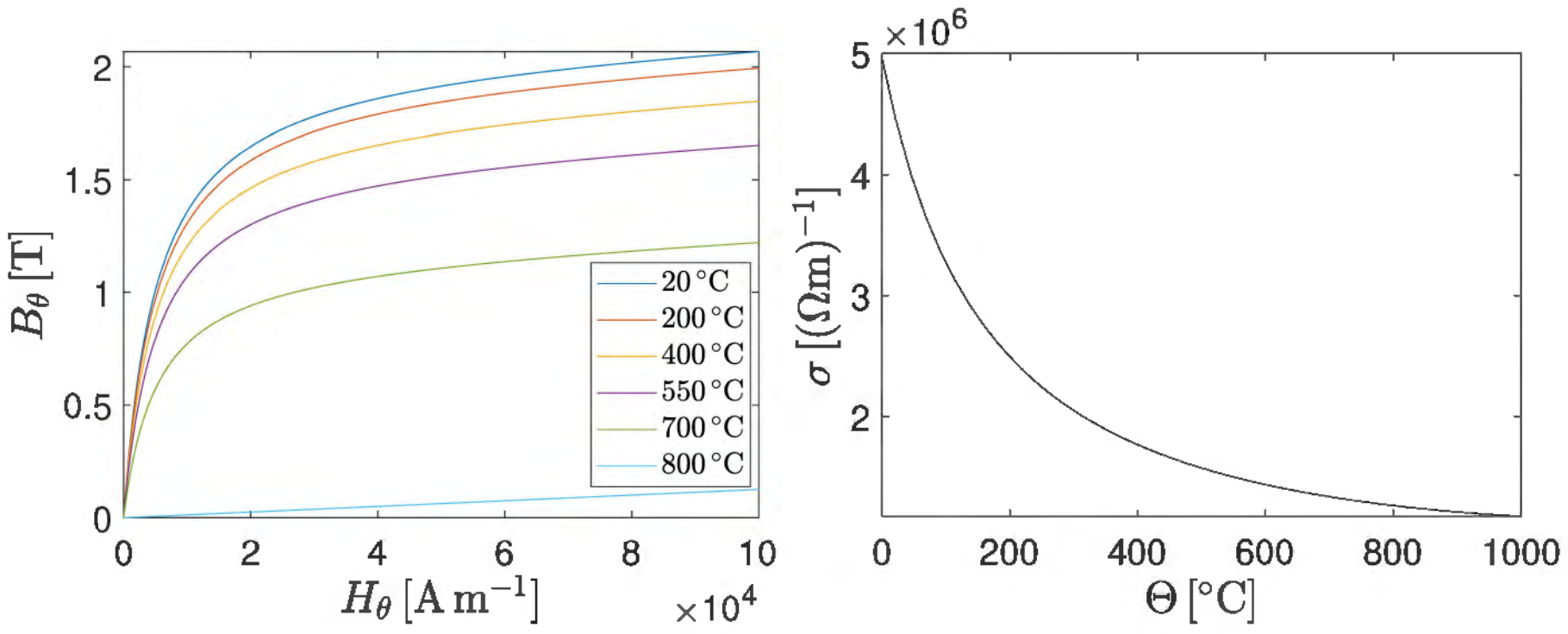

(a) Relationship between $H_\theta$ and $B_\theta$ at different temperatures (left) and electrical conductivity vs. temperature (right).

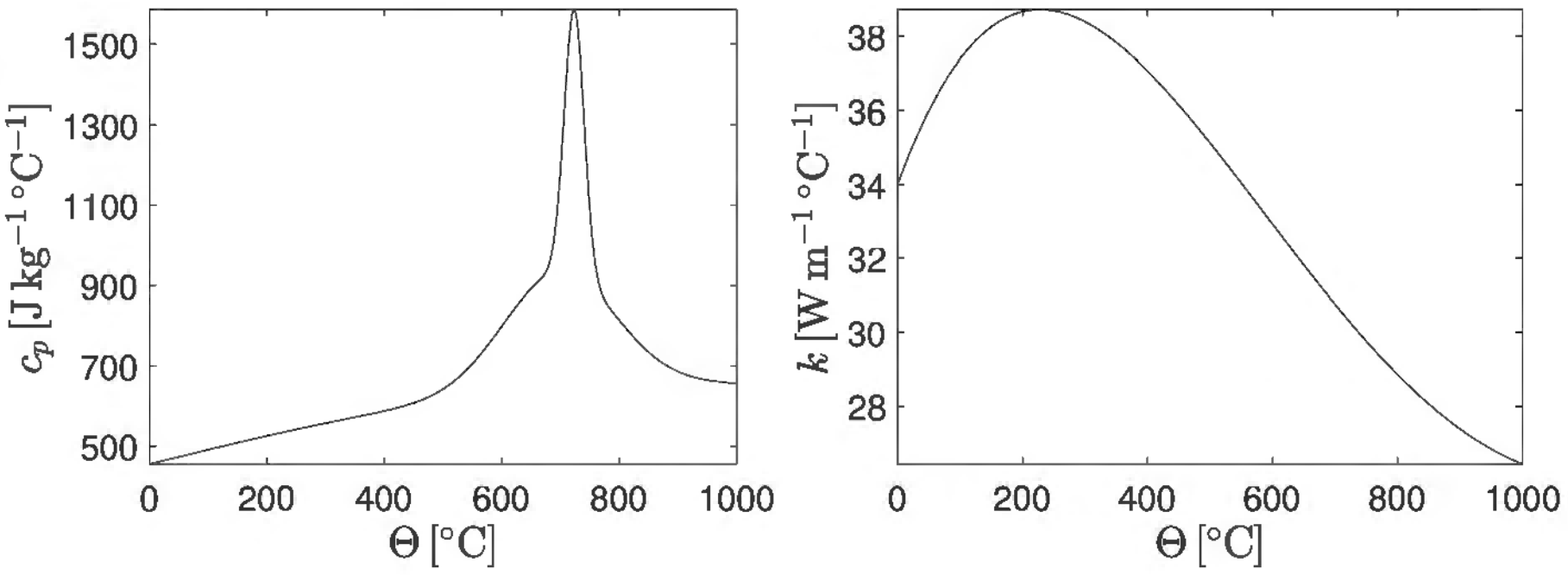

(b) Specific heat (left) and thermal conductivity (right) vs. temperature.

**Fig. 5**: Material properties.

problem was solved during 20 seconds of simulation and results are presented both for a time when the temperature is below the Curie temperature and for a later time when the Curie temperature has been reached in certain parts and the skin effect has therefore been attenuated.

Figure 6 and Figure 7 show the modulus of the field $\tilde{H}_\theta$ in Eulerian and Lagrangian approaches after 2 and 20 seconds, respectively.

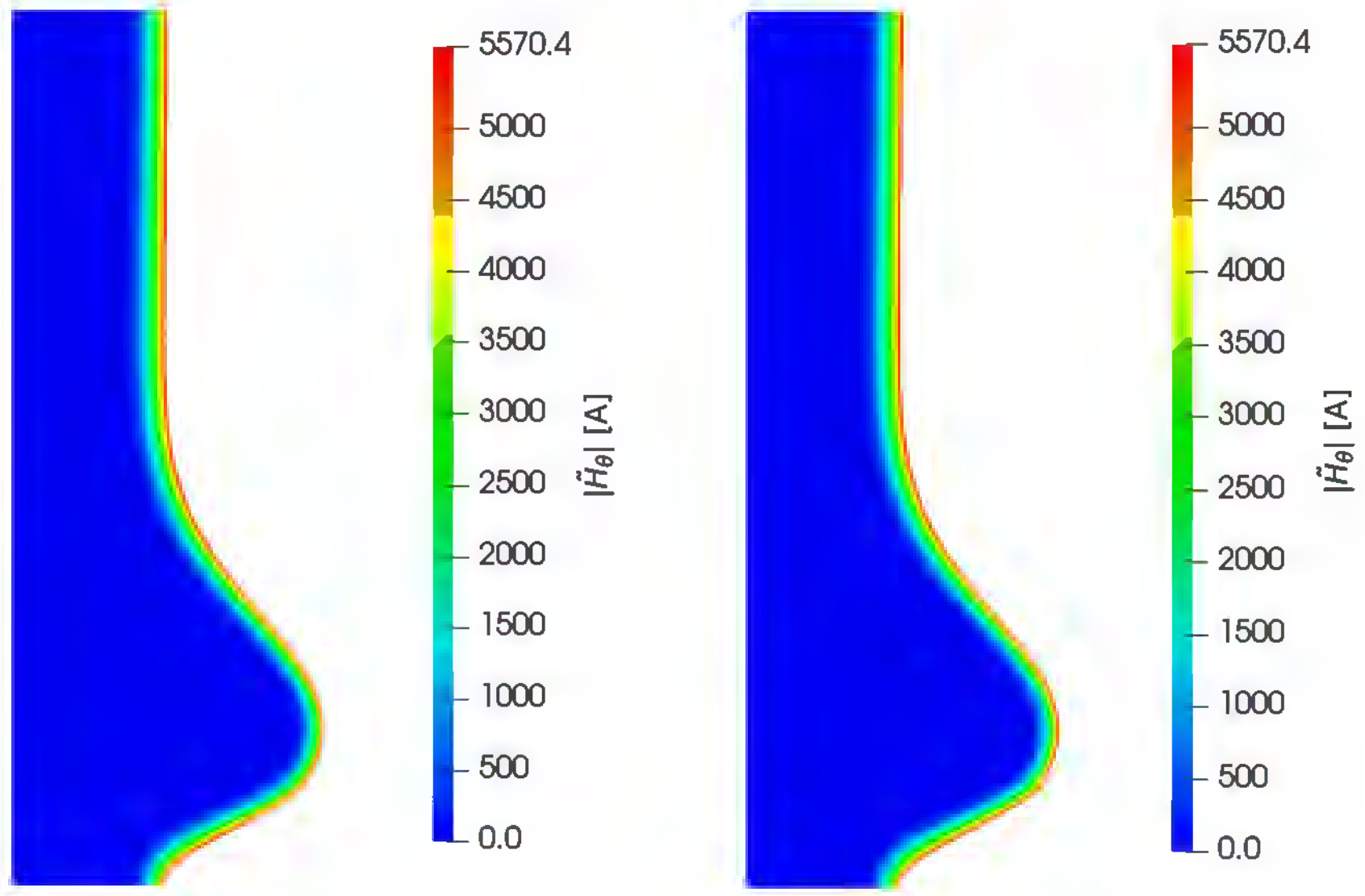

**Fig. 6**: Modulus of the field $\tilde{H}_\theta$ at $t = 2\,\mathrm{s}$: Eulerian approach (left) and Lagrangian approach (right).

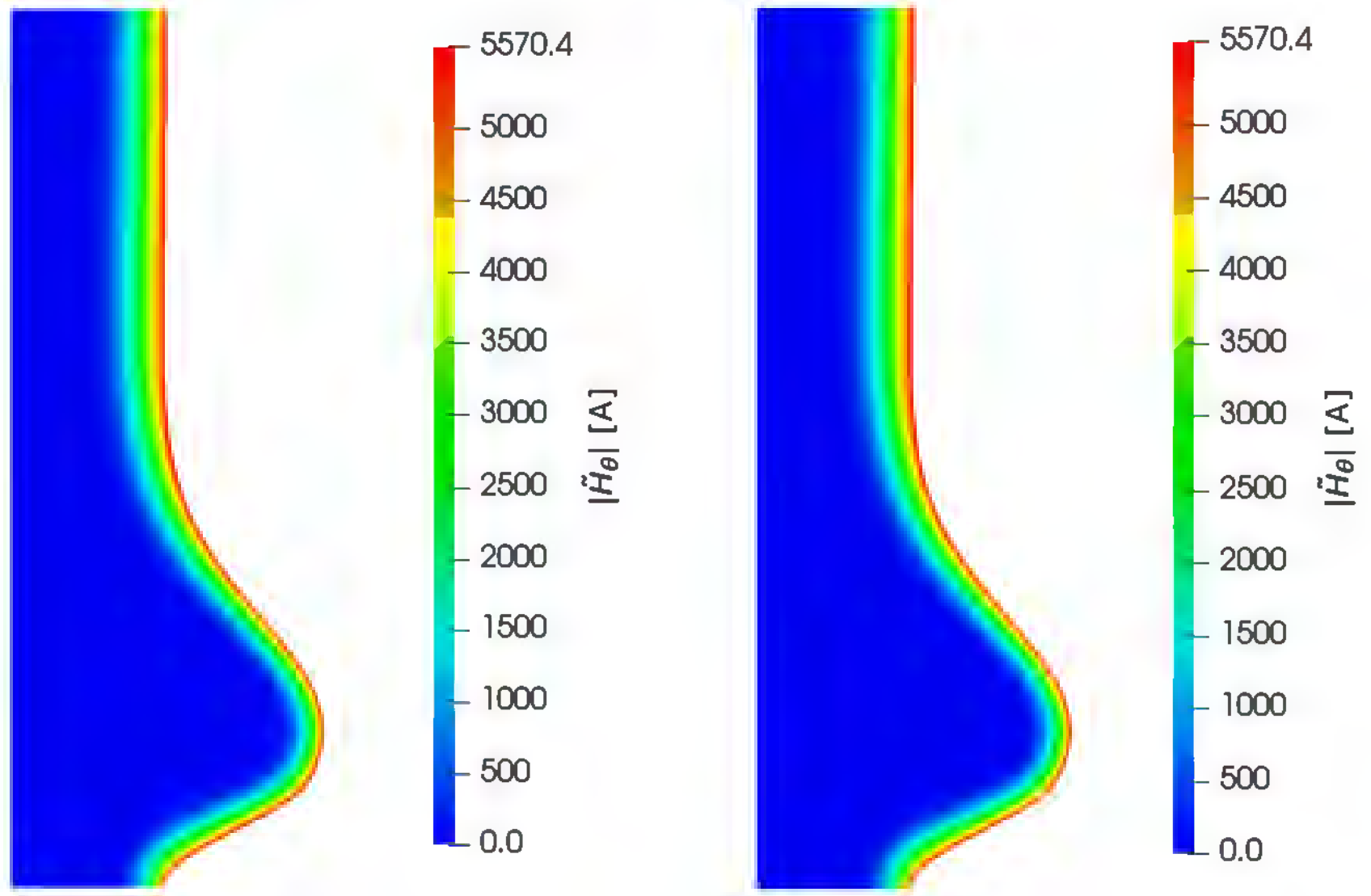

**Fig. 7**: Modulus of the field $\tilde{H}_\theta$ at $t = 20\,\mathrm{s}$: Eulerian approach (left) and Lagrangian approach (right).

On the other hand, the modulus of the electrical current density at different times is shown in Figure 8 and Figure 9.

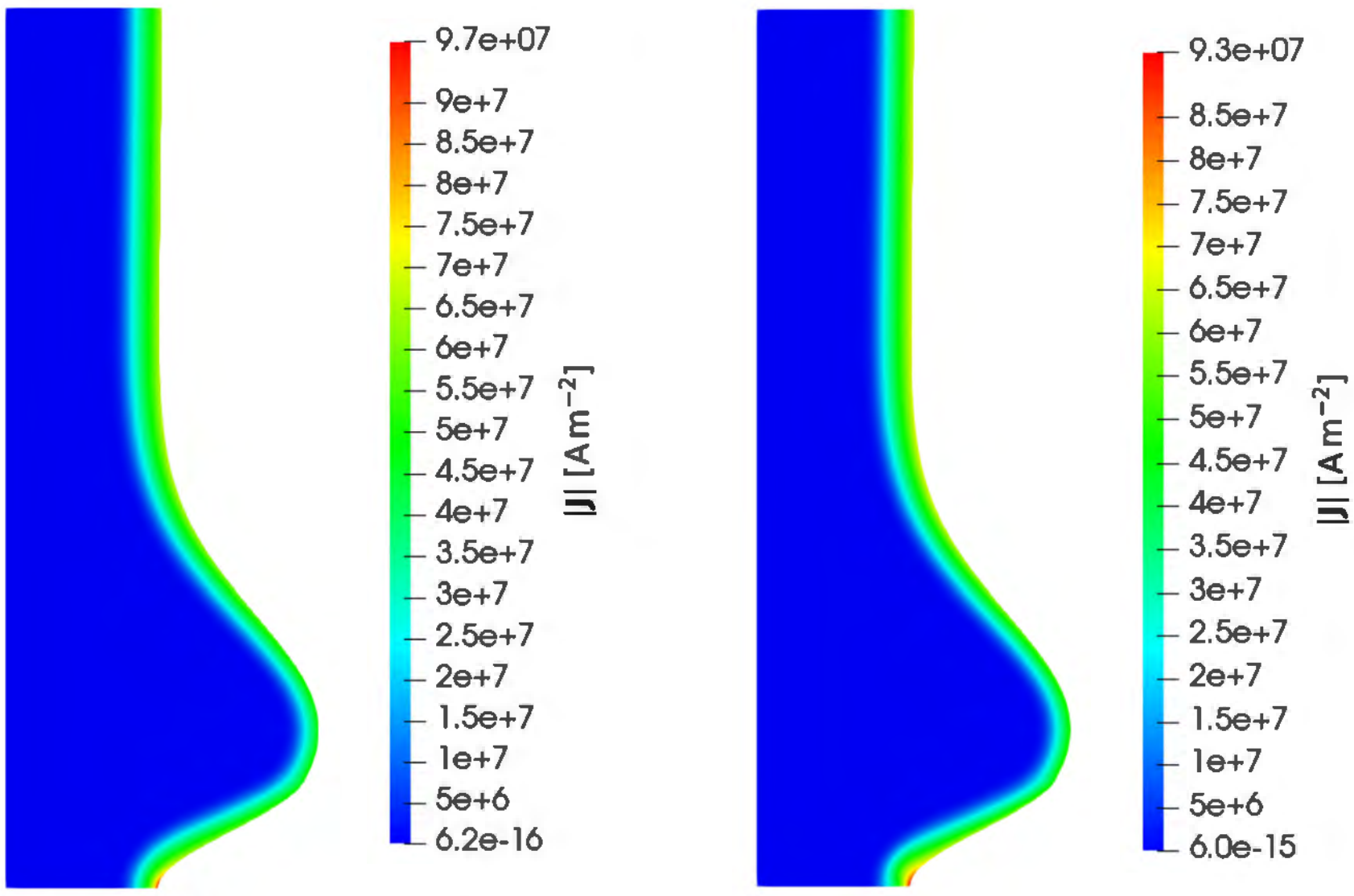

**Fig. 8**: Current density at $t = 2\,\mathrm{s}$: Eulerian approach (left) and Lagrangian approach (right).

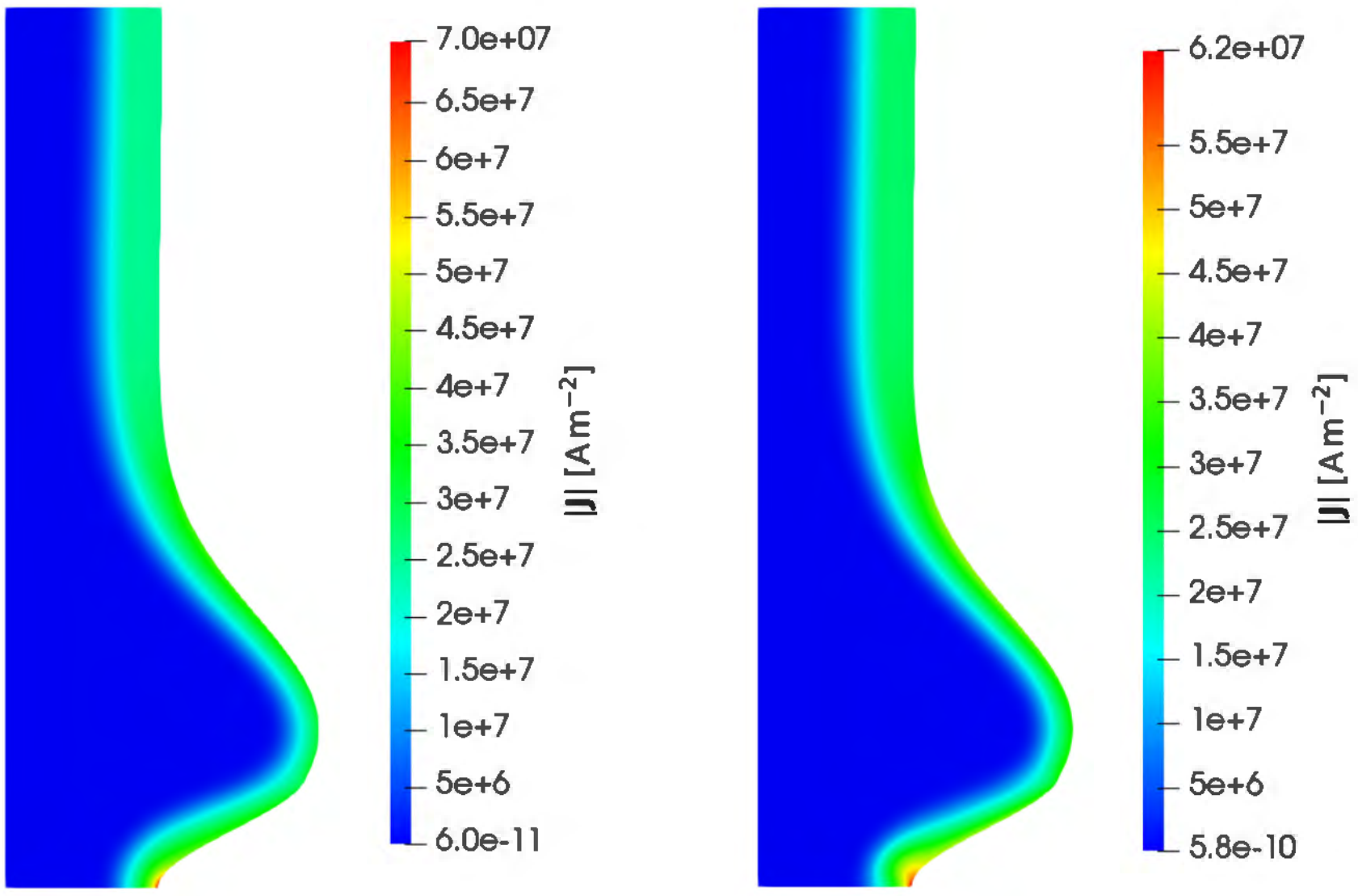

**Fig. 9**: Current density at $t = 20\,\mathrm{s}$: Eulerian approach (left) and Lagrangian approach (right).

The comparison of the voltage between the electric ports (top and bottom of the cylinder) is presented in Figure 10.

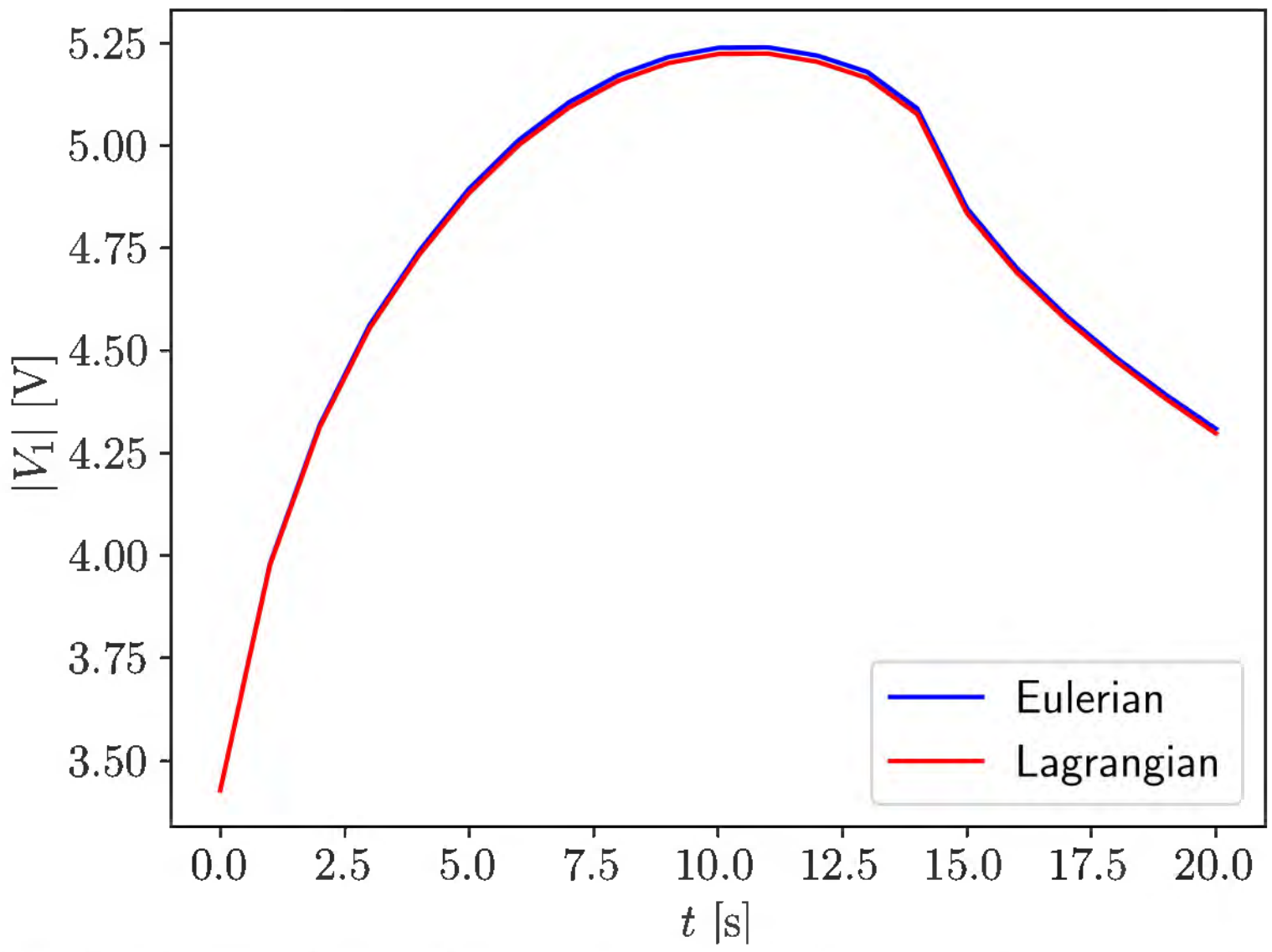

**Fig. 10**: Comparison of potential drop between the electric ports vs. time.

Finally, concerning the thermal part, the temperature field at different times is presented in Figure 11 and Figure 12.

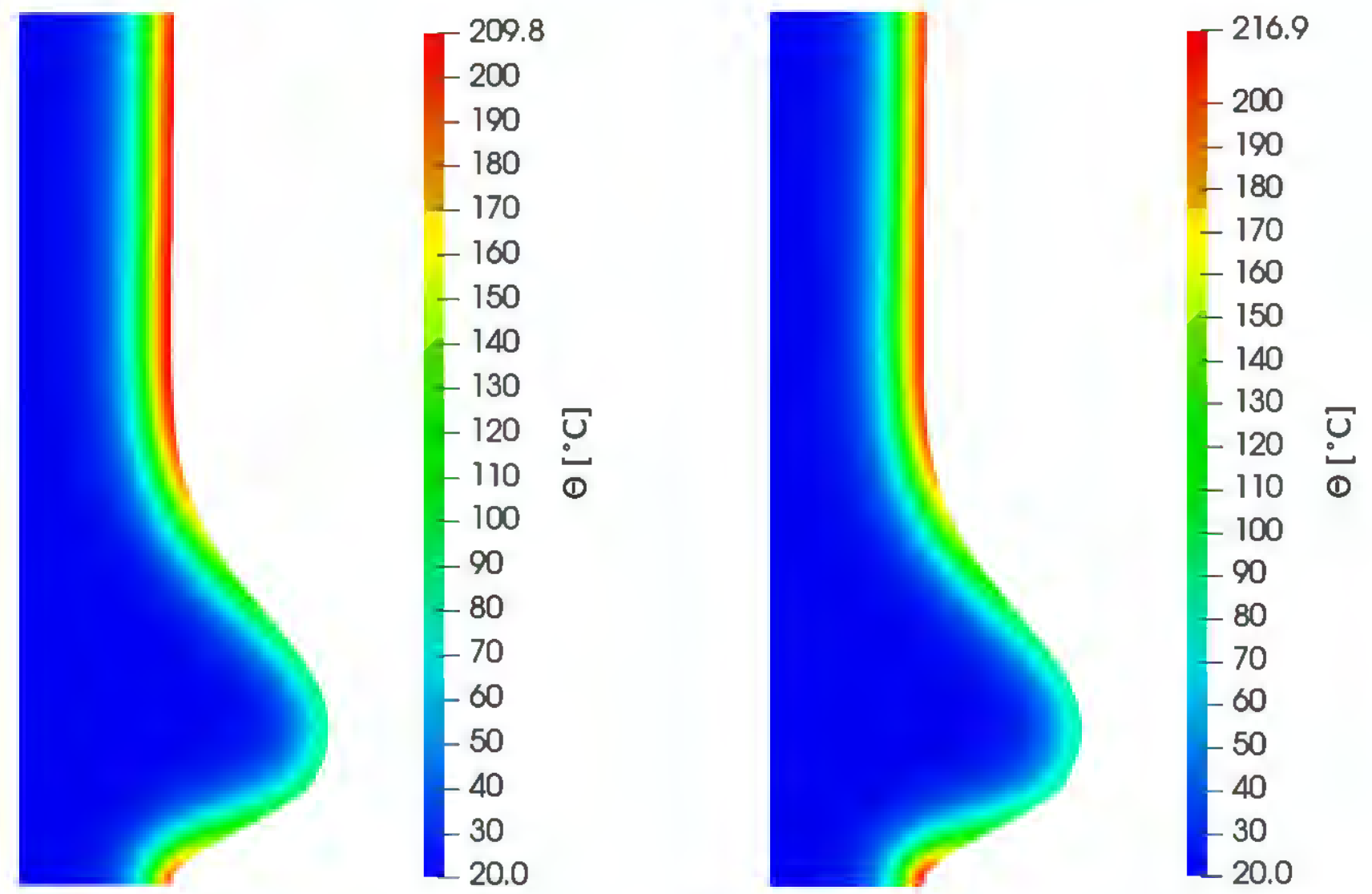

**Fig. 11**: Temperature field at $t = 2\,\mathrm{s}$: Eulerian approach (left) and Lagrangian approach (right).

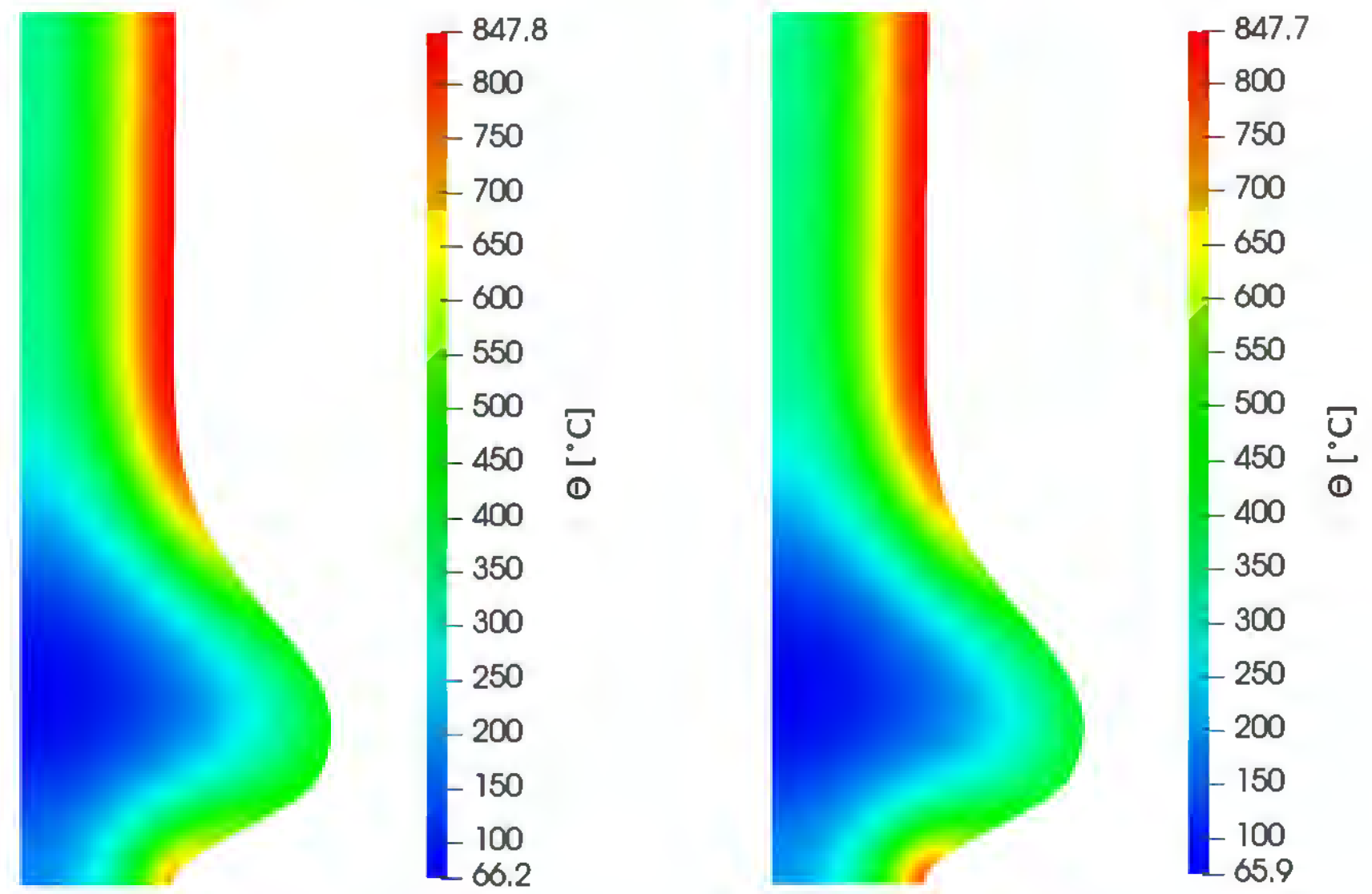

**Fig. 12**: Temperature field at $t = 20\,\mathrm{s}$: Eulerian approach (left) and Lagrangian approach (right).

Note that there is a good agreement between the results in Eulerian and Lagrangian coordinates.

# 6 Conclusions

A thermo-electromagnetic model for calculating the dissipated power and the temperature in cylindrical pieces undergoing large deformations has been introduced. A fully coupled problem is addressed, where a time-harmonic eddy current model with in-plane currents and electric ports is considered for the electromagnetic part, and a heat transfer transient model is used to describe the thermal part. Both models are approached from a Lagrangian point of view, with appropriate weak formulations that are later implemented in Python-FEniCS based on finite element methods. This approach is validated through numerical simulation by performing a suitable test also solved in Eulerian coordinates.

The forthcoming work, now in progress, is to couple the thermo-electromagnetic model to a non-linear mechanical model, which will allow us to consider large deformations and calculate the displacement field at each time step. This model will provide a powerful tool to simulate processes such as electric upsetting with an alternating source for cylindrical pieces. The methodology can be extended to deal with other forming processes involving large deformations, such as electromagnetic forming.

# Declarations

**Funding.** The research has been developed in collaboration with CIE Galfor through a project granted by the Centre for the Development of Industrial Technology (CDTI) and signed between the company CIE Galfor and Itmati (nowadays, integrated in CITMAga). This work has been partially supported by Xunta de Galicia under grant 2021 GRC GI-1563 ED431C 2021/15 and MCIN/AEI/10.13039/501100011033/ FEDER Una manera de hacer Europa, under the research project PID2021-122625OB-I00.

**Conflict of interest.** The authors declare that they have no conflict of interest.